\documentclass{amsart}
\usepackage{epsf}
\usepackage[all]{xy}
\usepackage{enumerate}
\usepackage{amssymb}
\usepackage{color}

\newcommand{\ie}{{\em i.e.\ }}

\newtheorem*{theorem}{Theorem}
\newtheorem*{lemma}{Lemma}
\newtheorem*{proposition}{Proposition}
\newtheorem*{corollary}{Corollary}

\theoremstyle{definition}

\newtheorem*{example}{Example}

\theoremstyle{remark}
\newtheorem*{remark}{\textbf{Remark}}

\numberwithin{equation}{section}

\newcommand{\internalcomment}[1]{}

\newcommand{\Z}{\mathbf{Z}}

\newcommand{\ko}{\: , \;}

\newcommand{\ol}{\overline}
\newcommand{\ul}{\underline}
\newcommand{\we}{\wedge}
\newcommand{\che}{\vee}
\renewcommand{\tilde}[1]{\widetilde{#1}}

\newcommand{\ra}{\rightarrow}

\newcommand{\arr}[1]{\stackrel{#1}{\rightarrow}}

\newcommand{\opname}[1]{\operatorname{\mathsf{#1}}}

\newcommand{\Mod}{\opname{Mod}\nolimits}
\newcommand{\Tria}{{\mbox{Tria}}}

\newcommand{\id}{\mathbf{1}}

\newcommand{\ten}{\otimes}

\newcommand{\cok}{\opname{cok}\nolimits}
\newcommand{\im}{\opname{im}\nolimits}
\renewcommand{\ker}{\opname{ker}\nolimits}

\newcommand{\colim}{\varinjlim}

\newcommand{\op}[1]{\opname{#1}\nolimits}

\renewcommand{\H}[1]{{H}^{#1}}

\newcommand{\cc}{{\mathcal C}}
\newcommand{\cd}{{\mathcal D}}
\newcommand{\ce}{{\mathcal E}}
\newcommand{\cf}{{\mathcal F}}

\newcommand{\ch}{{\mathcal H}}
\newcommand{\ci}{{\mathcal I}}
\newcommand{\cj}{{\mathcal J}}

\newcommand{\cn}{{\mathcal N}}

\newcommand{\cq}{{\mathcal Q}}

\newcommand{\cw}{{\mathcal W}}
\newcommand{\cx}{{\mathcal X}}

\newcommand{\eps}{\varepsilon}

\renewcommand{\phi}{\varphi}

\newcommand{\Hom}{\opname{Hom}}
\newcommand{\HOM}{\opname{Hom^\bullet}}

\newcommand{\cone}{\opname{Cone}\nolimits}

\newcommand{\centeps}[1]{\begin{array}{c} \epsfbox{#1} \end{array}}

\begin{document}
\title{The bar derived category of a curved dg algebra}

\author{Pedro Nicol\'{a}s}
\address{Departamento de Matem\'{a}ticas, Universidad de Murcia, Aptdo. 4021, 30100, Espinardo, Murcia, Espa\~na}

\email{pedronz@um.es}

%\author{Author Two}
%\address{Mathematical Research Section, School of Mathematical Sciences,
%Australian National University, Canberra ACT 2601, Australia}

%    General info
\date{February 7, 2008}



\begin{abstract}
Curved $A_{\infty}$-algebras appear in nature as deformations of dg
algebras. We develop the basic theory of curved $A_\infty$-algebras
and, in particular, curved dg algebras. We investigate 
their link with a suitable class of dg coalgebras via the bar 
construction and produce Quillen model structures on
their module categories. We define the analogue of
the relative derived category for a curved dg algebra.
\end{abstract}

\maketitle
%\tableofcontents

\section{Introduction}

\subsection{Motivation: deformations of dg categories}
The infinitesimal deformation theory of abelian categories and their
Hochschild cohomology have recently been established by Wendy T.
Lowen and Michel Van den Bergh \cite{LowenVandenBergh2004a,
LowenVandenBergh2004b}. This theory is motivated by non commutative
algebraic geometry and in particular the need to give a theoretical
framework for the \emph{ad hoc} arguments used in the construction of
important classes of non commutative projective varieties. The further
development of this theory requires a good understanding of the
deformation theory of `non commutative schemes'. These can be modeled
by differential graded (=dg) categories \cite{Drinfeld2004, Toen2006,
  Keller2006b, Keller2006c}, and this is the motivation for this
paper. The initial observation is that the full Hochschild complex of
a dg algebra does not parametrize the deformations in the category of
dg algebras but rather in the category of curved $A_\infty$-algebras,
which we also call $A_{[0,\infty[}$-algebras.  Key examples of these
are curved dg algebras, \ie graded algebras endowed with a degree one
derivation whose square is not necessarily zero but equals the
commutator with a given element of degree $2$. These algebras were
introduced by L.  E. Positsel'ski{\u\i} \cite{Positselskii1993}, who
showed that they occur in nature when one generalizes Koszul duality
to non homogeneous quadratic algebras.

$A_{[0,\infty[}$-algebras appear in the work of Ezra Getzler and
John D. S. Jones \cite{GetzlerJones1989} and the work of Ezra
Getzler, John D. S. Jones and Scott Petrack
\cite{GetzlerJonesPetrack1991} on $S^1$-equivariant differential
forms on the free loop space of a smooth manifold, in the work of Gunnar Fl\o ystad \cite{Floystad2006} on Koszul duality and in the work of Alberto S. Cattaneo and Giovanni Felder
\cite{CattaneoFelder2007} on the relative version of Kontsevich's
formality theorem. In Physics, curved dg algebras appear in the work
of Albert Schwarz \cite{Schwarz1999, Schwarz2003} as
noncommutative generalizations of  $Q$-manifolds, in the work of
Anton Kapustin and Yi Li \cite{KapustinLi2003a, KapustinLi2003b} on topological $D$-branes in Landau-Ginzburg
models, in the work of  Hiroshige Kajiura \cite{Kajiura2005} on
deformation of holomorphic line bundles over higher dimensional
complex tori, in the work of Xiang Tang \cite{Tang2006} on strict
quantization \dots

\subsection{Contents}
The aim of this paper is to develop the basic theory of curved
$A_\infty$-algebras and curved dg algebras, their module categories,
their link with a suitable class of dg coalgebras via the bar
construction, and to investigate the analogue of the relative derived
category for a curved dg algebra.

More precisely, we extend certain results on $A_\infty$-algebras and
their modules obtained in \cite{Lefevre03} in two directions:
\begin{itemize}
\item[-] instead of working over a field we work over an arbitrary commutative
ground ring;
\item[-] instead of considering $A_\infty$-algebras and dg algebras
we more generally consider $A_{[0,\infty[}$-algebras and curved dg
algebras.
\end{itemize}
Our results concern the bar/cobar adjunction for algebras and
modules and the existence of Quillen model structures on suitable
categories of modules.

In section \ref{Model structures on Frobenius categories}, we study
how to construct certain Quillen model structures (\emph{cf.}\ 
\cite{Hovey1999, Hirschhorn2003}) in Frobenius categories
(\emph{cf.}\  \cite{Keller1996}) inspired by the techniques of \cite[Section
2.3]{Hovey1999}. In particular, let
\[\xymatrix{\cc\ar@<1ex>[d]^{R} \\
\cd\ar@<1ex>[u]^L
}
\]
be an adjoint pair of functors and let $\eta:\id\ra RL$ and
$\delta:LR\ra\id$ be the  adjunction morphisms. Assume that $\cc$
and $\cd$ are Frobenius categories and $L$ and $R$ are exact
functors (and so they preserve injectives). Put $\cw_{\cc}$ (resp.
$\cw_{\cd}$) for the class of morphisms of $\cc$ (resp. $\cd$)
mapped to a stable isomorphism by $R$ (resp. $L$). We prove that if
$R(\delta_{M})$ and $L(\eta_{N})$ are stable isomorphisms, then 
\begin{enumerate}[a)]
\item There is a model structure in $\cc$ (resp. in $\cd$) having $\cw_{\cc}$ (resp. $\cw_{\cd}$) as the class of weak equivalences. In some sense, these are the minimal model structures such that the morphism $\delta_{M}:LRM\ra M$ is a cofibrant approximation and $\eta_{N}:N\ra RLN$ is a fibrant approximation.
\item The localizations $\cc[\cw^{-1}_{\cc}]$ and $\cd[\cw^{-1}_{\cd}]$ exist and they are triangulated quotients of the corresponding stable categories. The pair of adjoint functors $(L,R)$ induces mutually quasi-inverse triangulated equivalences
\[\xymatrix{\cc[\cw_{\cc}^{-1}]\ar@<1ex>[d] \\
\cd[\cw_{\cd}^{-1}]\ar@<1ex>[u]
}
\]
\item If both $\cc$ and $\cd$ have finite direct limits and colimits, then they are model categories endowed with those model structures, and the pair of adjoint functors $(L,R)$ is a Quillen equivalence.
\item If $\cc$ has kernels and $\cd$ has cokernels, then they are relevant derivable categories (\emph{cf.}\  \cite{Cisinski2003}) and the pair of adjoint functors $(L,R)$ induces mutually quasi-inverse equivalences between the corresponding derivators (\emph{cf.}\  \cite{Maltsiniotis2001, Maltsiniotis2005}).
\end{enumerate}

In sections \ref{A zero infinity algebras and their bar construction}, \ref{Bar cobar
adjunction for cdg algebras}, \ref{A zero infinity modules and their
bar construction} and \ref{Bar cobar adjunction for cdg modules} we
extend the bar/cobar formalism (\emph{cf.}\  \cite{Lefevre03} for this
formalism in the absence of the curvature, and
\cite{Positselskii1993, GetzlerJones1989, CattaneoFelder2007} for some of this formalism in the presence
of a non vanishing curvature) to $A_{[0,\infty[}$-algebras and their
modules, stressing the importance of Maurer-Cartan equations. It
turns out that these constructions somehow behave more naturally
when applied to a dg algebra not regarded as without curvature but
rather as a curved dg algebra with zero curvature. In section \ref{A
zero infinity algebras and their bar construction}, we give the
basic definitions concerning $A_{[0,\infty[}$-algebras, introduce
the cocomplete graded-augmented counital dg coalgebras and present
the bar construction of an $A_{[0,\infty[}$-algebra $A$ as a
representative of the functor which takes such a coalgebra $C$ to
the set of `twisting cochains' $\op{Tw}(C,A)$. In section \ref{Bar
cobar adjunction for cdg algebras}, we define the cobar construction
of a cocomplete graded-augmented counital dg coalgebra $C$ as a
corepresentative of the functor which takes a curved dg algebra $A$
to the set of `twisting cochains' $\op{Tw}(C,A)$. It is a left
adjoint to the restriction of the bar construction to the category
of curved dg algebras. In particular, for each
$A_{[0,\infty[}$-algebra $A$, we have a canonical curved dg algebra
$\Omega BA$ and a natural morphism of $A_{[0,\infty[}$-algebras
$A\ra\Omega B A$ universal among the morphisms of
$A_{[0,\infty[}$-algebras from $A$ to a curved dg algebra. In
section \ref{A zero infinity modules and their bar construction}, we
present the basic definitions concerning $A_{[0,\infty[}$-modules,
define the `linearized' Maurer-Cartan equation of an
$A_{[0,\infty[}$-module and construct (in great generality) the bar
construction of an $A_{[0,\infty[}$-module $M$ as a representative
of a functor which takes certain comodules $N$ to the `linearized'
Maurer-Cartan equation of an $A_{[0,\infty[}$-module defined from
$M$ and $N$. In section \ref{Bar cobar adjunction for cdg modules},
we define the cobar construction of a counital dg comodule $N$ as a
corepresentative of the functor which takes a unital curved dg
module $M$ to the `linearized' Maurer-Cartan equation of a curved dg
module defined from $M$ and $N$. This is a left adjoint of the
restriction of the bar construction to the category of unital curved
dg modules. From this adjunction we introduce in section \ref{The
bar derived category} the \emph{bar derived category} of a unital
curved dg algebra as a certain homotopy category, by using some
results of section \ref{Model structures on Frobenius categories}.
In section \ref{The various derived categories of a dg algebra}, we
prove that the bar derived category of a curved dg algebra with zero
curvature is the relative derived category (\emph{cf.}\ \cite{Keller1998a})
of the underlying dg algebra. In particular, this gives a model
structure for the relative derived category, proves that the bar
derived category of a dg algebra over a field coincides with the
classical derived category and allows us a better comprehension of
Kenji Lef{\`e}vre-Hasegawa's theorem \cite[Th\'{e}or\`{e}me
2.2.2.2]{Lefevre03}. Some results of \cite{Lefevre03} suggest that
the spirit of $A_{\infty}$-theory is to replace quasi-isomorphisms
by homotopy equivalences, up to increasing the amount of morphisms
and/or objects. In section 10, we show that the $A_{\infty}$-theory
over an arbitrary commutative ring still allows us to describe the
relative derived category of an augmented dg algebra $A$ as the
category of unital dg $A$-modules, with (strictly unital) morphisms
of $A_{\infty}$-modules, up to (strictly unital) homotopy
equivalences of $A_{\infty}$-modules.

\subsection{Acknowledgments}
I thank Bernhard Keller and Wendy Tor Lowen for helpful and enthusiastic discussions on the topic of this article. Also, I acknowledge the hospitality of the Institut de Math\'ematiques de Jussieu.

\section{Notation}

Unless otherwise stated, $k$ will be a commutative (associative, unital) ring. Also, `graded' will always mean `$\Z$-graded'. If $V$ is a graded $k$-module, \ie
\[V = \bigoplus_{p\in\Z} V^p \ko
\]
we denote by $SV$ or $V[1]$ the graded $k$-module with $(SV)^p=V^{p+1}$ for all $p\in\Z$. We call $SV$ the \emph{suspension} or the
\emph{shift} of $V$. The shift extends to an automorphism of the category of graded $k$-modules, with inverse denoted by $S^{-1}$. Notice that, given a graded $k$-module $V$, we have two homogeneous morphisms $s:V\ra SV$ and $w:V\ra S^{-1}V$, of degree $-1$ and $1$ respectively.

If $f: U \to U'$ and $g: V \to V'$ are homogeneous morphisms between graded $k$-modules, their \emph{tensor product}
\[ f \ten g : U \ten V \to U'\ten V'
\]
is defined using the \emph{Koszul sign rule}: We have
\[ (f \ten g)(u \ten v)= (-1)^{|g|\cdot|u|} f(u) \ten g(v)
\]
for all homogeneous elements $v\in V$ and $w\in W$, where $|g|$ and $|v|$ are the degrees of $g$ and $v$, respectively.

An \emph{exact category} in the sense of Quillen \cite{Quillen73} is an additive category $\cc$ endowed with a distinguished class $\ce$ of sequences  
\[X\arr{i}Y\arr{p}Z,
\]
closed under isomorphisms such that $(i,p)$ is an \emph{exact pair}, \ie $i$ is the kernel of $p$ and $p$ is the cokernel of $i$. Following \cite{GabrielRoiter}, the morphisms $p$ are called \emph{deflations}, the morphisms $i$ \emph{inflations} and the pairs $(i, p)$ \emph{conflations}. The class of conflations have to satisfy the following axioms:
\begin{enumerate}
\item[Ex0\ )] The identity morphism of the zero object is a deflation.
\item[Ex1\ )] The composition of two deflations is a deflation.
\item[Ex1')] The composition of two inflations is an inflation.
\item[Ex2\ )] Deflations admit and are stable under base change.
\item[Ex2')] Inflations admit and are stable under cobase change.
\end{enumerate}
As shown by B. Keller \cite{Keller1990}, these axioms are equivalent to Quillen's and they imply that if $\cc$ is
small, then there is a fully faithful functor from $\cc$ into an ambient abelian category $\cc'$ whose image
is an additive subcategory closed under extensions and such that a sequence of $\cc$ is a
conflation if and only if its image is a short exact sequence of $\cc'$. Conversely, one easily checks that
an extension closed full additive subcategory $\cc$ of an abelian category $\cc'$ endowed with all exact pairs which induce short exact sequences in $\cc'$ is always exact.

A \emph{Frobenius category} is an exact category $\cc$ with enough $\ce$-injectives and enough $\ce$-projectives and where the class of $\ce$-projectives coincides with the class of $\ce$-injectives. In this case, the \emph{stable category} $\ul{\cc}$ obtained by dividing $\cc$ by the ideal of morphisms factoring
through an $\ce$-projective-injective carries a canonical structure of triangulated category, \emph{cf.}\ \cite{Heller68, Happel87, KellerVossieck87, GelfandManin}. We write $\ol{f}$ for the image in $\ul{\cc}$ of a morphism $f$ of $\cc$. The \emph{suspension or shift functor} $S$ of $\ul{\cc}$ is obtained by choosing a conflation
\[X\arr{i_{X}}IX\arr{p_{X}}SX
\]
for each object $X$, where $IX$ is required to be $\ce$-injective. Each triangle is isomorphic to a standard triangle $(\ol{i},\ol{p},\ol{e})$ obtained by embedding a conflation $(i,p)$ into a commutative diagram
\[\xymatrix{X\ar[r]^{i}\ar[d]^{\id} & Y\ar[r]^{p}\ar[d] & Z\ar[d]^{e} \\
X\ar[r]^{i_{X}} & IX\ar[r]^{p_{X}} & SX
}
\]

For the notation concerning \emph{model categories} we refer to \cite{Hovey1999, Hirschhorn2003}. However, our notion of \emph{model structure} is weaker than that of \cite{Hovey1999} since we do not impose functorial factorizations. We will say that a \emph{model structure} on a category $\cc$ is the data of three classes of morphisms, $\cw\ko \cc of$ and $\cf ib$, the \emph{weak equivalences}, \emph{cofibrations} and \emph{fibrations} respectively, satisfying the following axioms:
\begin{enumerate}[1)]
\item (2-out-of-3) If $f$ and $g$ are morphisms of $\cc$ such that $gf$ is defined and two of $f, g$ and $gf$ are weak equivalences, then so is the third.
\item (Retracts) The three classes of morphisms are closed under retracts (in the category of morphisms of $\cc$).
\item (Lifting) Define a map to be a \emph{trivial cofibration} if it is both a cofibration and a weak equivalence. Similarly, define a map to be a \emph{trivial fibration} if it is both a fibration and a weak equivalence. Then trivial cofibrations have the left lifting property with respect to fibrations, and cofibrations have the left lifting property with repect to trivial fibrations.
\item (Factorization) For any morphism $f$ of $\cc$ there exist factorizations
\[\xymatrix{X\ar[rr]^{f}\ar[dr]_{\alpha} && Y \\
& Y'\ar[ur]_{\beta} &
}\ko
\xymatrix{X\ar[rr]^f\ar[dr]_{\gamma} && Y \\
& X'\ar[ur]_{\delta} &
}
\]
where $\alpha$ is a cofibration, $\beta$ is a trivial fibration, $\gamma$ is a trivial cofibration and $\delta$ is a fibration.
\end{enumerate}
We will say that a \emph{model category} is a category with finite direct and inverse limits endowed with a model structure. 

A category has a model structure if and only if it satisfies axioms M2), M5) and M6) of \cite{Quillen1967}. Hence, if a category has a model structure and finite direct and inverse limits, it is what D. Quillen called a \emph{closed model category} in \cite{Quillen1967}. 

For the more general notion of \emph{derivable category} we refer to \cite{Cisinski2003}. In short, a derivable category is a category $\cc$ endowed with three classes of morphisms $\cw\ko \cc of$ and $\cf ib$ satisfying some familiar axioms which ensure that its \emph{homotopy category} $\op{Ho}\cc:=\cc[\cw^{-1}]$, \ie the localization of $\cc$ with respect to the class $\cw$, is still understandable, that is to say, it is a category whose morphisms can be calculated by homotopy relations and calculus of fractions. A \emph{relevant derivable category} (\emph{cf.}\ subsection 5.1. of \cite{Cisinski2003}) is a derivable category satisfying a certain `lifting property'. One of the good properties of relevant derivable categories is that the morphisms in the homotopy categories can be calculated just by using homotopy relations \cite[Proposition 5.11]{Cisinski2003}, and so the homotopy categories have small $\Hom$-sets.

Notice also that if a category has a model structure, an initial object, a final object, pullbacks and pushouts, then it is a relevant derivable category.

%\begin{proof}
%Let us see that if $\cc$ is such a category, then it satisfies the axioms D0)-D4) of left derivable category and their duals.

%D0) $\cc$ has a final object by hypothesis. By the retract argument, we have that $\cf ib=T\cc of$-inj, and so $\cf ib$ is closed under compositions. Then, to see that the class of fibrant objects is closed under isomorphism it suffices to see that isomorphisms are fibrations. But, since $T\cf ib=\cc of$-inj, we can easily prove that isomorphisms are trivial fibrations. On the other hand, if $1$ is the final object, since $\id_{1}$ is trivial fibration, then $1$ is fibrant.

%D1) It was already proved in D0) that every isomorphism is a weak equivalence. The second part of D1) is the 2-out-of-3 axiom.

%D2) It was already proved in D0) that every isomorphism is a trivial fibration and that the class of fibrations is closed under compositions. The last part of axiom D2) is also easy to prove since, by hypothesis, $\cc$ has pullbacks, and it is well known that the class of fibrations is closed under base change since $\cf ib=T\cc of$-inj.

%D3) Since $T\cf ib=\cc of$-inj, we also have that the class of trivial fibrations is closed under base change.

%D4) is part of the Factorization Axiom of the definition of `model structure'.

%Dually, $\cc$ satisfies D0$^{op}$)-D4$^{op})$. Finally, thanks to the Lifting Axiom of the definition of `model structure', $\cc$ is a relevant derivable category.
%\end{proof}
\bigskip

\section{Model structures on Frobenius categories}\label{Model structures on Frobenius categories}

\subsection{A recognition criterion in Frobenius categories}\label{A recognition criterion in Frobenius categories}
Recall that if $\cx$ is a class of morphisms of a category $\cc$, then $\cx$-inj is the class formed by those morphisms of $\cc$ having the right lifting property with respect to every morphism in $\cx$, $\cx$-proj is the class formed by those morphisms of $\cc$ having the left lifting property with respect to every morphism in $\cx$ and $\cx$-cof$=(\cx$-inj$)$-proj.

The proof of \cite[Theorem 2.1.19]{Hovey1999} gives us a kind of `recognition criterion' to detect model structures. Indeed, consider three classes of morphisms, $\cw\ko \ci$ and $\cj$ in a category $\cc$. Then, there is a model structure on $\cc$ with $\ci$-cof as the class of cofibrations, $\cj$-cof as the class of trivial cofibrations, and $\cw$ as the class of weak equivalences if
\begin{enumerate}[1)]
\item $\cw$ has the $2$-out-of-$3$ property and is closed under retracts.
\item Any morphism of $\cc$ factors as a morphism in $\ci$-cof followed by a morphism in $\ci$-inj.
\item Any morphism of $\cc$ factors as a morphism in $\cj$-cof followed by a morphism in $\cj$-inj.
\item $\cj$-cof $\subseteq\cw\cap\ci$-cof.
\item $\ci$-inj $\subseteq\cw\cap\cj$-inj.
\item Either $\cw\cap\ci$-cof $\subseteq\cj$-cof or $\cw\cap\cj$-inj $\subseteq\ci$-inj.
\end{enumerate}

The following result is a generalization of the techniques of \cite[Section 2.3]{Hovey1999}. Thanks to the recognition criterion, it tells us that if in a Frobenius category we find a class $\cq$ of objects (closed under shifts) such that the class $\ci$ of morphisms $i_{Q}:Q\ra IQ$ with $Q$ in $\cq$, and the class $\cj$ of morphisms $0\ra IQ$ with $Q$ in $\cq$, satisfy the conditions 2) and 3) above, then we have a model structure.

\begin{theorem}
Let $\cc$ be a k-linear Frobenius category and let $\cq$ be a class
of objects of $\cc$ closed under shifts. Consider the classes $\ci$
formed by the morphisms $i_{Q}:Q\ra IQ$ where $Q$ runs through the
class $\cq$, $\cj$ formed by the morphisms $0\ra IQ$ where $Q$ runs
through the class $\cq$ and $\cw$ of morphisms $f$ such that
$\ul{\cc}(Q,\ol{f})$ is an isomorphism of $k$-modules for all the
objects $Q$ in $\cq$. If $\ci$ (resp. $\cj$) allows factorizations
of any morphism as a map in $\ci$-cof (resp. $\cj$-cof) followed by
a map in $\ci$-inj (resp. $\cj$-inj), then there is a model
structure on $\cc$ with $\ci$-cof as the class of cofibrations,
$\cj$-cof as the class of trivial cofibrations and $\cw$ as the
class of weak equivalences. Moreover, with respect to this model
structure we have that
\begin{enumerate}[1)]
\item every object is fibrant,
\item every cofibration is an inflation,
\item every inflation with cokernel in $\cq$ is a cofibration,
\item two morphisms are left homotopic if and only if their difference factors through an injective.
\end{enumerate}
\end{theorem}
\begin{proof}
Define $\cf ib:=\cj$-inj and $\cc of:=\ci$-cof. The elements of $\cf ib$ will be called \emph{fibrations}, and the terms of $\cc of$ will be called \emph{cofibrations}.

$\bullet$ Let us prove that if $p:X\ra Y$ is a morphism of $\ci$-inj, then the induced map $p^{\che}:\cc(Q,X)\ra\cc(Q,Y)$ is surjective for all $Q$ in $\cq$. Given a morphism $y:Q\ra Y$, we consider the (solid) commutative square with its (dotted) lifting
\[\xymatrix{S^{-1}Q\ar[r]^{0}\ar[d]_{i_{S^{-1}Q}} & X\ar[d]^{p} \\
IS^{-1}Q\ar[r]_{\ \ \ yp_{S^{-1}Q}}\ar@{.>}[ur]_{z} & Y
}
\]
Since $zi_{S^{-1}Q}=0$, there exists a unique morphism $z':Q\ra X$ such that $z'p_{S^{-1}Q}=z$, and so from $pz=yp_{S^{-1}Q}$ we deduce $pz'=y$.

$\bullet$ Let us prove the inclusion $\ci\text{-inj}\subseteq\cj\text{-inj}$. Indeed, given a (commutative) diagram
\[\xymatrix{0\ar[r]\ar[d] & X\ar[d]^{p} \\
IQ\ar[r]_{y} & Y
}
\]
with $p$ in $\ci$-inj, since $p^{\che}:\cc(Q,X)\ra\cc(Q,Y)$ is surjective, there exists a morphism $x:Q\ra X$ such that $px=yi_{Q}$. Now, the following (solid) commutative diagram has a (dotted) lifting
\[\xymatrix{Q\ar[r]^{x}\ar[d]_{i_{Q}}& X\ar[d]^{p} \\
IQ\ar[r]_{y}\ar@{.>}[ur]_{z} & Y
}
\]

$\bullet$ Hence $\ci\text{-inj}\subseteq\cj$-inj, and therefore $\cj\text{-cof}\subseteq\ci$-cof.

$\bullet$ Let us prove that if $p:X\ra Y$ is a morphism in $\ci$-inj, then the induced map $\ol{p}^{\che}:\ul{\cc}(Q,X)\ra\ul{\cc}(Q,Y)$ is bijective for all $Q$ in $\cq$. It is surjective because at the level of Frobenius categories it was surjective. Assume that $x:Q\ra X$ is a morphism such that the composition $px: Q\ra Y$ factors through an injective. Then it is forced to factor also through $i_{Q}$, and so it is of the form $px=yi_{Q}$. Then, the lifting in the following commutative diagram tells us that $x$ factors through an injective:
\[\xymatrix{Q\ar[r]^{x}\ar[d]_{i_{Q}} & X\ar[d]^{p} \\
IQ\ar[r]_{y}\ar@{.>}[ur]_{z} & Y
}
\]

$\bullet$ Hence $\ci\text{-inj}\subseteq\cw\cap(\cj\text{-inj})$.

$\bullet$ Let us prove that if $p: X\ra Y$ is a morphism in $\cw\cap(\cj\text{-inj})$, then $\ul{\cc}(Q,\ker p)=0$ for all $Q$ in $\cq$. Observe that $p\in\cj$-inj means that $p^{\che}:\cc(IQ,X)\ra\cc(IQ,Y)$ is surjective for all $Q$ in $\cq$. Hence, by diagram chasing in the following commutative diagram with exact
columns and exact rows
\[\xymatrix{\cc(IQ,X)\ar[r]^{p^{\che}}\ar[d]^{(i_{Q})^{\we}} & \cc(IQ,Y)\ar[r]\ar[d]^{(i_{Q})^{\we}} & 0 \\
\cc(Q,X)\ar[r]^{p^{\che}}\ar[d] & \cc(Q,Y)\ar[d] & \\
\ul{\cc}(Q,X)\ar[r]^{\ol{p}^{\che}}\ar[d] & \ul{\cc}(Q,Y)\ar[r]\ar[d] & 0 \\
0 & 0 &
}
\]
we deduce that $p^{\che}:\cc(Q,X)\ra\cc(Q,Y)$ is surjective for all $Q$ in $\cq$. Now, we apply the Snake Lemma to the commutative diagram with exact rows
\[\xymatrix{0\ar[r] & \cc(IQ,\ker p)\ar[r]^{(p^k)^{\che}}\ar[d]^{(i_{Q})^{\we}} & \cc(IQ,X)\ar[r]^{p^{\che}}\ar[d]^{(i_{Q})^{\we}} & \cc(IQ,Y)\ar[r]\ar[d]^{(i_{Q})^{\we}} & 0 \\
0\ar[r] & \cc(Q,\ker p)\ar[r]^{(p^k)^{\che}} & \cc(Q,X)\ar[r]^{p^{\che}} & \cc(Q,Y)\ar[r] & 0
}
\]
and we deduce the following long exact sequence
\[\dots\ra\cc(SQ,X)\arr{p^{\che}}\cc(SQ,Y)\arr{\delta}\ul{\cc}(Q,\ker p)\arr{(\ol{p^k})^{\che}}\ul{\cc}(Q,X)\arr{(\ol{p})^{\che}}\ul{\cc}(Q,Y)\ra 0
\]
Since $p^{\che}$ is surjective (as proved before), then $\delta=0$. Since $(\ol{p})^{\che}$ is injective, then $(\ol{p^k})^{\che}=0$. Therefore, $\ul{\cc}(Q,\ker p)=0$.

$\bullet$ Let us prove that if $p\in\cw\cap(\cj\text{-inj})$, then $p\in\ci$-inj. Consider the following commutative diagram
\[\xymatrix{Q\ar[r]^{x}\ar[d]_{i_{Q}} & X\ar[d]^{p} \\
IQ\ar[r]_{y} & Y
}
\]
Since $p$ belongs to $\cj$-inj, there exists a morphism $w:IQ\ra X$ such that $pw=y$. Now $p(wi_{Q}-x)=0$, and so $wi_{Q}-x=p^ku$ for some morphism $u:Q\ra\ker p$. But since
$\ul{\cc}(Q,\ker p)=0$, there exists $v:IQ\ra \ker p$ such that $vi_{Q}=u$. Thus, $(w-p^kv)i_{Q}=x$ and $p(w-p^kv)=y$.

$\bullet$ We have $\ci\text{-inj}=\cw\cap(\cj\text{-inj})$

$\bullet$ Let us prove that every morphism in $\ci\text{-cof}$ is an inflation. Given $f:X\ra Y$ in $\ci\text{-cof}$, we consider a lifting in the following diagram
\[\xymatrix{ X\ar[r]^{i_{X}}\ar[d]_{f} & IX\ar[d] \\
Y\ar[r]\ar@{.>}[ur]_{g} & 0
}
\]
This proves that $f$ is a monomorphism in $\cc$ (and in the ambient abelian category \cite[Proposition A.2]{Keller1990}, since $i_{X}$ is a monomorphism in that abelian category), and we get the following commutative diagram
\[\xymatrix{X\ar[r]^{f}\ar[d]_{\id_{X}} & Y\ar[r]^{f^c}\ar[d]^g & \cok f\ar[d]^{e} \\
X\ar[r]_{i_{X}} & IX\ar[r]_{p_{X}} & SX
}
\]
Since the right square is cocartesian, then $f^c$ is a deflation. Since $f$ is a mono in the ambient abelian category, then $f$ is the kernel of $f^c$, and so
$f$ is an inflation.

$\bullet$ Let us prove that if $f$ is a morphism in $\cj\text{-cof}$, then its cokernel $\cok f$ is projective relative to fibrations. Indeed, let $p:M\ra N$ be a fibration and consider a morphism
$g:\cok f\ra N$. We form the following (solid) commutative diagram and consider its (dotted) lifting
\[\xymatrix{ X\ar[r]^0\ar[d]_{f} & M\ar[d]^p \\
Y\ar[r]_{gf^c}\ar@{.>}[ur]^{h} & N
}
\]
Since $hf=0$, then there exists a unique $h':\cok f\ra M$ such that $h'f^c=h$. From $ph=gf^c$ we deduce $ph'f^c=gf^c$, and so $ph'=g$.

$\bullet$ In particular, since $\cf ib=\cj$-inj contains all the deflations, the cokernel of a morphism in $\cj$-cof is projective(-injective).

$\bullet$ Given a morphism $f$ in $\cj$-cof, since $\cj\text{-cof}\subseteq\ci$-cof, we have that $f$ is an inflation. Hence the conflation $X\arr{f}Y\arr{f^c}\cok f$ yields a triangle
\[X\arr{\ol{f}}Y\arr{\ol{f^c}}\cok f\ra SX
\]
But $\cok f$ is projective, and so $\ol{f}$ is an isomorphism. In particular, $f\in\cw$.

$\bullet$ Then $\cj\text{-cof}\subseteq\cw\cap(\ci\text{-cof})$.

$\bullet$ We have proved that $\ci\text{-inj}=\cw\cap(\cj\text{-inj})$ and $\cj\text{-cof}\subseteq\cw\cap(\ci\text{-cof})$. Therefore, the hypotheses of the recognition criterion are satisfied.

$\bullet$ Now, notice that for an arbitrary object $M$ the morphisms
\[\xymatrix{M\oplus M\ar[rrr]^{\scriptsize{\left[\begin{array}{cc}i_{M}&0 \\ 0 & i_{M} \\ \id_{M} & \id_{M}\end{array}\right]}\ \ \ } &&& \ IM\oplus IM\oplus M\ \ \ar[rrr]^{\ \ \ \scriptsize{\left[\begin{array}{ccc}0&0&\id_{M}\end{array}\right]}} &&& M
}
\]
form a factorization of $\scriptsize{\left[\begin{array}{cc}\id_{M}&\id_{M}\end{array}\right]}$ ending in a weak equivalence, and so it is a cylinder object. Then, if two morphisms $f,g:M\ra N$ are left homotopic then their difference factors through an injective. On the other hand, if $f-g=hi_{M}$ for some $h$, then by considering $\scriptsize{\left[\begin{array}{ccc}h&0&g\end{array}\right]}$ we prove that $f$ and $g$ are left homotopic.

That every object is fibrant is obvious, and that every cofibration is
an inflation has been seen when proving the existence of the model
structure.

Now, let $f:X\ra Y$ be an inflation with cokernel $Q$ in $\cq$, and let $g:M\ra N$ be a trivial fibration. We fit any commutative square of the form
\[\xymatrix{X\ar[r]\ar[d]^f & M\ar[d]^g \\
Y\ar[r] & N
}
\]
in a commutative diagram
\[\xymatrix{S^{-1}Q\ar[r]\ar[d]_{h} & X\ar[r]\ar[d]^f & M\ar[d]^g \\
PQ\ar[r]\ar[d] & Y\ar[r]\ar[d] & N \\
Q\ar@{=}[r] & Q &
}
\]
where the top left square is bicartesian. Since $g$ has the right lifting property with respect to $h$, then the same holds with respect to $f$. This proves that
$f\in\ci$-cof.
\end{proof}

\begin{lemma}
Let $\cc$ be a category with initial and final object. Assume that $\cc$ is endowed with a model structure $(\cw,\cf ib,\cc of)$ such that the class of cofibrations admit cobase changes. If $\cc$ has pullbacks, then it is a relevant derivable category.
\end{lemma}
\begin{proof}
By using the retract argument and  \cite[Lemma 7.2.11]{Hirschhorn2003}, one shows easily that a category with a final object and pullbacks endowed with a model structure is a left derivable category. Let us see that $\cc$ is also a right derivable category. Indeed, the only axioms which are not trivially satisfied are D2$^{op}$) and D3$^{op}$), but they hold since cofibrations admit cobase changes. 
\end{proof}

\begin{corollary}
Let $\cc$ be a $k$-linear Frobenius category and let $\cq$ be a class of objects of $\cc$ closed under shifts. Assume that $\cq$ is such that induces a model structure as in the theorem above. If $\cc$ has kernels, then it is a relevant derivable category. In particular, its homotopy category $\op{Ho}\cc:=\cc[\cw^{-1}]$ has small $\Hom$-sets.
\end{corollary}
\begin{proof}
Since $\cc$ has finite products and kernels, then it has pullbacks. On the other hand, inflations admit cobase changes and cofibrations are inflations. Therefore, we can apply the lemma above.
\end{proof}

In what follows, we present two situations in which a class $\cq$ of objects allows good factorizations.

\subsection{Small object argument in Frobenius categories}\label{Small object argument in Frobenius categories}
If $\cq$ is a set (resp. a class) such that $\ci$ (resp. and $\cj$)
allows the (generalized) small object argument \cite{Hirschhorn2003}
\cite{Chorny06}, then we easily get a model structure.

The following result is a generalization of \cite[Theorem 2.3.11]{Hovey1999}.

\begin{corollary}
Let $\cc$ be a $k$-linear Frobenius category with small colimits.
Let $\cq$ be a set of objects of $\cc$ closed under shifts. Define
$\ci$ to be the class of morphisms $i_{Q}:Q\ra IQ$ where $Q$ runs
through the set $\cq$, $\cj$ the class of morphisms $0\ra IQ$ where
$Q$ runs through $\cq$ and $\cw$ the class of morphisms $f$ such
that $\ul{\cc}(Q,\ol{f})$ is an isomorphism for all $Q$ in $\cq$.
Assume that the objects of $\cq$ are small relative to $\ci$-cell
(\emph{cf.}\ \cite[Definition 10.4.1.]{Hirschhorn2003}). Then the classes
above define a model structure on $\cc$, and so, if $\cc$ has kernels, it is a relevant derivable
category. Moreover, its homotopy category $\op{Ho}\cc$ is the
triangle quotient of $\ul{\cc}$ by the full triangulated subcategory
whose objects are the $M$ such that $\ul{\cc}(Q,M)=0$ for all $Q$ in
$\cq$, and is triangle equivalent to $\Tria(\cq)$, the smallest
full triangulated subcategory of $\ul{\cc}$ containing $\cq$ and
closed under coproducts.
\end{corollary}
\begin{proof}
  The hypotheses 2) and 3) of the recognition criterion are
  satisfied thanks to the small object argument \cite[Proposition
  10.5.16, Proposition 10.5.10]{Hirschhorn2003}. Hence, $\cc$ has a
  model structure. To prove the last statement, consider the full
  triangulated subcategory $\cn$ of $\ul{\cc}$ whose objects are the
  $M$ such that $\ul{\cc}(Q,M)=0$ for all $Q$ in $\cq$.
  Let us show that $(\Tria(\cq),\cn)$ is a t-structure on $\ul{\cc}$. For this, it
  is enough to see that for each object $M$, there exists a triangle
\[
M\ra M'\ra M''\ra SM
\]
in $\ul{\cc}$ with $M'\in\cn$ and $M''\in \Tria(\cq)$. Now, given $M$, thanks to the factorization associated to $\ci$, we have a
relative $\ci$-cell complex $f: M\ra M'$ such that $M'\ra 0$ is in $\ci$-inj, \ie $M'\in\cn$. Since every relative $\ci$-cell complex is
a cofibration \cite[Proposition 10.5.10]{Hirschhorn2003} and
cofibrations are inflations, we have a triangle $M\arr{\ol{f}} M'\ra
M''\ra SM$ coming from a conflation $M\arr{f} M'\ra M''$. 
Since $f$ is a relative $\ci$-cell complex, its cokernel $M''$ is an $\ci$-cell complex and, by the lemma below, it belongs to
$\Tria(\cq)$. 

Therefore, $(\Tria(\cq), \cn)$ is a t-structure
and we have a series of equivalences
\[
\Tria(\cq)\simeq \ul{\cc}/\cn\simeq \cc[\cw^{-1}]=\op{Ho}\cc.
\]
\end{proof}

\begin{lemma}
Under the hypotheses of the corollary, each $\ci$-cell complex belongs to $\Tria(\cq)$.
\end{lemma}

\begin{proof}

\noindent \emph{First step: Let $\lambda$ be an ordinal. If we have
a direct system of conflations
\[
\eps_\alpha: 0 \to X_\alpha \to Y_\alpha \to Z_\alpha \to 0 \ko
\alpha<\lambda \ko
\]
such that the structure morphisms $Z_\alpha \to Z_{\beta}$ are
inflations for all $\alpha<\beta<\lambda$, then the colimit of the
system is a conflation.} Indeed, it suffices to check that for each
injective $I$, the sequence of abelian groups
\[
0 \to \cc(\colim Z_\alpha, I) \to \cc(\colim Y_\alpha, I) \to
\cc(\colim X_\alpha, I) \to 0
\]
is exact. This follows from the Mittag-Leffler criterion \cite[0$_{\text{III}}$, 13.1]{EGAIIIa} since the
maps
\[
\cc(Z_\beta, I) \to \cc(Z_\alpha, I)
\]
are surjective for all $\alpha<\beta<\lambda$.

\noindent \emph{Second step: If we have an acyclic complex of $\cc$
\[
\ldots \to X^p \to X^{p+1} \to \ldots \to X^0 \to Y \to 0\ko
\]
then $Y$ belongs to the smallest triangulated subcategory of
$\ul{\cc}$ containing the $X^p$ and stable under countable
coproducts.} Indeed, by lemma~6.1 of \cite{Keller1990}, for each
complex $K$ over $\cc$, there is a triangle
\[
\mathbf{a} K \to K \to \mathbf{i} K \to S \mathbf{a} K
\]
of $\ch(\cc)$ such that $\mathbf{i}K$ has injective components and
$\mathbf{a}K$ is the colimit (in the category of complexes) of a
countable sequence of componentwise split monomorphisms of acyclic
complexes. The functor $\mathbf{i}K$ is the left adjoint of the
inclusion into $\ch(\cc)$ of the full subcategory of complexes with
injective components. Thus, it commutes with coproducts. The
composed functor
\[
F: \ch(\cc) \to \ul{\cc} \ko K \mapsto Z^0(\mathbf{a} K)
\]
is a triangle functor which commutes with coproducts and extends the
projection $\cc\to \ul{\cc}$ from $\cc$ to $\ch(\cc)$. Moreover, it
vanishes on acyclic complexes. Thus, it maps the truncated complex
\[
X'=(\ldots \to X^p \to X^{p+1} \to \ldots \to X^0 \to 0 \to \ldots)
\]
to an object isomorphic to $Y$ in $\ul{\cc}$. Since $F$ commutes
with coproducts, it suffices to show that the complex $X'$ is in the
smallest triangulated subcategory of $\ch(\cc)$ containing the $X^p$
and stable under countable coproducts. This holds thanks to Milnor's
triangle (\emph{cf.}\ \cite{Keller1998b, Milnor1962})
\[
\coprod X^{\geq p} \to \coprod X^{\geq q} \to X' \to S\coprod
X^{\geq p} \ko
\]
where $X^{\geq p}$ is the subcomplex
\[
0 \to X^p \to X^{p+1} \to \ldots \to X^0 \to 0
\]
and the leftmost morphism has the components
\[
\xymatrix{X^{\geq p} \ar[r]^-{[\id, -i]^t} & X^{\geq p}\oplus
X^{\geq p+1} \ar[r] & \coprod X^{\geq q} } \ko
\]
where $i$ is the inclusion $X^{\geq p} \to X^{\geq p+1}$.

\noindent \emph{Third step: The claim.} Let $X$ be an $\ci$-cell complex. Then there is an ordinal $\lambda$ and a direct system
$X_\alpha$, $\alpha\leq\lambda$, such that we have $X=X_\lambda$ and
\begin{itemize}
\item[-] $X_0=0$,
\item[-] for all $\alpha<\lambda$, the morphism $X_\alpha\to
X_{\alpha+1}$ is an inflation with cokernel in $\cq$,
\item[-] for all limit ordinals $\beta\leq \lambda$, we have
$X_\beta = \colim_{\alpha<\beta} X_\alpha$.
\end{itemize}
We will show by induction on $\beta\leq \lambda$ that $X_\beta$ belongs to
$\Tria(\cq)$. This is clear for $X_0$. Moreover, if $X_\alpha$ is in
$\Tria(\cq)$, so is $X_{\alpha+1}$. So let us assume that $\beta$ is
a limit ordinal and $X_\alpha$ belongs to $\Tria(\cq)$ for each
$\alpha<\beta$. We wish to show that $X_\beta$ belongs to
$\Tria(\cq)$. Let $\cc^\beta$ be the category of functors $\beta \to
\cc$. The evaluation
\[
\cc^\beta \to \cc \ko Y \to Y_\alpha
\]
admits a left adjoint denoted by $Z\mapsto Z\ten \alpha$. For each
$Y\in \cc^\beta$, the morphism
\[
\coprod_{\alpha<\beta} Y_\alpha\ten\alpha \to Y
\]
is a pointwise split epimorphism. By splicing exact sequences of the form
\[
0 \to Y' \to \coprod_{\alpha<\beta} Y_\alpha\ten\alpha \to Y \to 0
\]
we construct a complex
\[
\ldots \to X^p \to X^{p+1} \to \ldots \to X^0 \to X \to 0
\]
which is ayclic for the pointwise split exact structure on
$\cc^\beta$ and such that each $X^p$ is a coproduct of objects
$Y\ten \alpha$, $Y\in \Tria(\cq)$, $\alpha<\beta$. By the first
step, the colimit $C$ of the above complex is still acyclic.
Moreover, the components of $C$ are coproducts of objects
\[
\colim_\gamma (Y\ten \alpha)(\gamma) = Y
\]
%In fact, $\colim:\cc^{\beta}\ra\cc$ is a left adjoint and so it preserves coproducts. Hence
%\[\colim\left(\coprod_{\alpha<\beta}X_{\alpha}\otimes\alpha\right)=\coprod_{\alpha<\beta}\left(\colim X_{\alpha}\otimes\alpha\right)=\coprod_{\alpha<\beta}X_{\alpha}.
%\]
belonging to $\Tria(\cq)$. Thus, each component of $C$ belongs to
$\Tria(\cq)$. Now, the claim follows from the second step.
\end{proof}

Notice that if in the corollary above, $\cc$ has finite inverse
(resp. small) limits, then it is a (resp. cofibrantly generated)
model category.

The fact that the pair $(\Tria(\cq),\cn)$ in the proof of Corollary \ref{Small object argument in Frobenius categories} is a
t-structure in $\ul{\cc}$ is a generalization of \cite[Proposition
4.5]{AlonsoJeremiasSouto2000}. With the same techniques, we can also
show that if $\cq$ is a set of objects closed under non-negative
shifts in a $k$-linear Frobenius category $\cc$ such that the
associated set $\ci$, formed by the maps $i_{Q}:Q\ra IQ$ where $Q$
runs through $\cq$, is small relative to $\ci$-cell, then the
smallest full suspended subcategory of $\ul{\cc}$ containing $\cq$
and closed under arbitrary coproducts is an aisle in $\ul{\cc}$.
This is a generalization of \cite[Proposition
3.2]{AlonsoJeremiasSouto2003}.

\begin{example}
Let $A$ be a unital dg $k$-algebra. The category of unital dg right $A$-modules, $\cc A$, has a structure of $k$-linear Frobenius category since it is the category of $0$-cocycles of a certain exact dg category, \emph{cf.}\ \cite{Keller1999, Keller2006b}. The conflations are those short exact sequences which split in the category of graded $A$-modules, and an object $M$ of $\cc A$ is injective-projective if and only if its identity morphism $\id_{M}$ is \emph{null-homotopic}, \ie $\id_{M}=d_{M}h+hd_{M}$ for some morphism $h$ of graded $A$-modules homogeneous of degree $-1$. The corresponding stable category is the \emph{category of unital dg $A$-modules up to homotopy}, $\ch A$. If we take $\cq$ to be the set formed by all the modules of the form $A[n]\ko n\in\Z$, then the smallness condition is satisfied (\emph{cf.}\ \cite[Example 2.1.6]{Hovey1999}) and the model structure leads to the derived category $\op{Ho}\cc A\cong\cd A$.
\end{example}

\begin{example}
  Let $A$ be a unital curved dg $k$-algebra with curvature $c\in A^2$
  (\emph{cf.}\ section \ref{A zero infinity algebras and their bar
    construction}). Let $\cc A$ be the category of unital curved dg
  $k$-modules (\emph{cf.}\ section \ref{Bar cobar adjunction for cdg
    modules}). One could define the graded $k$-module
\[\H \bullet(M):=\ker d_{M}^{\bullet}/(\im d_{M}^{\bullet-1}\cap\ker d^{\bullet}_{M})
\]
to be the \emph{naive cohomology} of a curved dg module $M$. Notice that $A$ is not a curved dg module with its regular structure, but
\[\tilde{A}:=A/cA
\]
is a curved dg right $A$-module with the natural multiplication by scalars. If we apply Corollary \ref{Small object argument in Frobenius categories} to the class $\cq$ formed by the objects $\tilde{A}[n]\ko n\in\Z$, then $\cc A$ becomes a model category whose weak equivalences are those morphisms inducing isomorphisms in naive cohomology. Its homotopy category would be a \emph{naive derived category}. Ezra Getzler and John D. S. Jones define in \cite{GetzlerJones1989} a unital curved dg algebra to be \emph{standard} if the curvature is in the center of the algebra. For these standard algebras they define a cohomology in the usual way. It is straightforward to prove that a morphism between standard unital curved dg algebras which induces an isomorphism at the level of that cohomology induces a triangle equivalence between the corresponding naive derived categories.
\end{example}

\begin{remark}
If the curvature of a unital curved dg algebra vanishes, then its naive derived category is precisely the classical derived category of the underlying dg algebra, and its \emph{bar derived category} (\emph{cf.}\ subsection \ref{Definition of the bar derived category}) is precisely the relative derived category of the underlying dg algebra. This seems to suggest that for an arbitrary unital curved dg algebra there would exist a fully faithful triangulated functor from its naive derived category to its bar derived category. However, the presence of a non-vanishing curvature makes things considerably less clear.
\end{remark}

\subsection{Factorizations provided by adjunctions}\label{Factorizations provided by adjunctions}
This section is motivated by the need of understanding the Quillen equivalence of Kenji Lef{\`e}vre-Hasegawa's theorem \cite[Th\'{e}or\`{e}me 2.2.2.2]{Lefevre03}, \emph{via} the bar/cobar constructions of section \ref{Bar cobar adjunction for cdg modules}.

Let $\cc$ and $\cd$ be $k$-linear Frobenius categories. We will prove that, under mild assumptions, an adjunction of exact functors 
\[\xymatrix{\cc\ar@<1ex>[d]^{R} \\
\cd\ar@<1ex>[u]^{L}
}
\]
(notice that this already implies that both $L$ and $R$ preserve injective objects) always gives rise to a class of objects $\cq$ of $\cc$ whose
associated classes of morphisms $\ci$ and $\cj$ of subsection \ref{A recognition criterion in Frobenius categories} allow the
factorizations required in conditions 2) and 3) of the recognition criterion. Dually, we also get a model structure in $\cd$
and, moreover, the localizations with respect to the corresponding classes of weak equivalences, $\cc[\cw^{-1}_{\cc}]$ and $\cd[\cw^{-1}_{\cd}]$, are triangle equivalent triangulated categories with small $\Hom$-sets. 

Let $\eta_{N}:N\ra RLN$ be the unit of the adjunction,
$\delta_{M}:LRM\ra M$ the counit of the adjunction,
\[\tau_{N,M}:\cc(LN,M)\arr{\sim}\cd(N,RM)
\]
the adjunction isomorphism.

\begin{theorem}\label{Factorizations provided by adjunctions}
If $R(\delta_{M})$ has injective kernel for each object $M$ of $\cc$, then $\cc$ admits a model structure in which
          \begin{enumerate}[a)]
          \item The cofibrations are the inflations with cokernel of the form $LN$ and their retracts.
          \item The fibrations are the morphisms with the right lifting property with respect to $0\ra PLN$ for every object $N$ of $\cd$.
          \item The class of weak equivalences, $\cw_{\cc}$, is formed by the morphisms $f$ such that $R(f)$ is a stable isomorphism.
          \end{enumerate}
         With respect to this model structure every object is fibrant and an object is cofibrant if and only if it is a direct summand of some $LN$. Also, we have the following
         description of the trivial cofibrations and the trivial fibrations:
         \begin{enumerate}[d)]
         \item The trivial cofibrations are the inflations with cokernel of the form $PLN$ and their retracts.
         \end{enumerate}
         \begin{enumerate}[e)]
         \item The trivial fibrations are the morphisms with the right lifting property with respect to the morphisms $S^{-1}LN\ra PLN$ for every object $N$ of $\cd$.
         \end{enumerate}
         Moreover, the localization of $\cc$ with respect to the weak equivalences is the triangle quotient of $\ul{\cc}$ by the full triangulated subcategory whose objects are the $M$ such that $\ul{\cc}(LN,M)=0$ for all $N$ in $\cd$, and is triangle equivalent to $\Tria(\{LN\}_{N\in\cd})$, the smallest full triangulated subcategory of $\ul{\cc}$ containing $LN\ko N\in\cd$ and closed under coproducts. The objects of this subcategory are precisely the direct summands, in $\ul{\cc}$, of the images $LN$ of objects $N$ of $\cd$ under $L$.
\end{theorem}
\begin{proof}
We apply Theorem \ref{A recognition criterion in Frobenius categories} to the class $\cq$ formed by the objects $LN$, where $N$ belongs to $\cd$. It is closed under shifts since $L$ commutes with the shift.

$\bullet$ Let $g:A\ra B$ be an inflation with cokernel of the form
$LN$. Let us show that $g$ is in $\ci$-cof. Indeed, if $f$ is a morphism of
$\ci$-inj, we can fit a commutative square of the form
\[\xymatrix{A\ar[r]\ar[d]^g & C\ar[d]^f \\
B\ar[r] & D
}
\]
into the commutative diagram
\[\xymatrix{ S^{-1}LN\ar[r]\ar[d]^h & A\ar[r]\ar[d]^g & C\ar[d]^f \\
PLN\ar[r]\ar[d] & B\ar[r]\ar[d] & D \\
LN\ar@{=}[r] & LN &
}
\]
Since $h$ belongs to $\ci$-cof and the left top square is bicartesian, then $g$ is in $\ci$-cof.

$\bullet$ Since $R(\delta_{M})$ is a retraction with injective kernel, we have that $\delta_{M}$ is in $\ci$-inj.

$\bullet$ \emph{Factorization associated to $\ci$:}
\[\xymatrix{ & X\ar[rr]^{f}\ar@{=}[dl]\ar[dr]_{\scriptsize{\left[\begin{array}{c}\iota_{X}\\ f\end{array}\right]}} & & Y  \\
X\ar[dr]_{\beta} & & IX\oplus Y\ar[ur]_{\scriptsize{\left[\begin{array}{cc}0 & \id\end{array}\right]}}\ar[dr]_{\scriptsize{\left[\begin{array}{cc}h & -g\end{array}\right]}} & \\
& E\ar[ur]^{\alpha}\ar[dr] & & C \\
&& LRC\ar[ur]_{\delta_{C}} &
}
\]
It has been constructed by forming a conflation over $\scriptsize{\left[\begin{array}{c}\iota_{X}\\ f\end{array}\right]}$ thanks to the diagram
\[\xymatrix{X\ar[r]^{\iota_{X}}\ar[d]^{f} & IX\ar[r]\ar[d]^{h} & SX\ar@{=}[d] \\
Y\ar[r]^{g} & C\ar[r] & SX
}
\]
made \emph{via} the pushout $C$, and then by forming the pullback $E$. We know already that $\beta$ is in $\ci$-cof and that
$\scriptsize{\left[\begin{array}{cc}0&\id\end{array}\right]}$ belongs to $\ci$-inj. Since $\delta_{C}$ is in $\ci$-inj,
then the composition $\scriptsize{\left[\begin{array}{cc}0&\id\end{array}\right]}\alpha$ is in $\ci$-inj.

$\bullet$ \emph{Factorization associated to $\cj$:}
\[\xymatrix{ X\ar[rr]^{f}\ar[dr]_{\scriptsize{\left[\begin{array}{c}0\\ \id\end{array}\right]}} && Y \\
& PLRY\oplus X\ar[ur]_{\ \ \scriptsize{\left[\begin{array}{cc}p_{Y}\delta_{PY}& f\end{array}\right]}} &
}
\]
formed by using the composition
\[\xymatrix{PLRY\oplus X=LRPY\oplus X\ar[rr]^{\ \ \ \ \ \ \ \ \ \ \scriptsize{\left[\begin{array}{cc}\delta_{PY} & 0\\ 0&\id\end{array}\right]}}
&&PY\oplus X\ar[rr]^{\scriptsize{\left[\begin{array}{cc}p_{Y}&f\end{array}\right]}} && Y
}
\]
It is clear that $\scriptsize{\left[\begin{array}{c}0\\ \id\end{array}\right]}$ belongs to $\cj$-cof and that
$\scriptsize{\left[\begin{array}{cc}\delta_{PY}& 0\\ 0 &\id\end{array}\right]}$ belongs to $\cj$-inj (since $\delta_{PY}$ is in $\ci$-inj$\subseteq\cj$-inj). Also, since
$\scriptsize{\left[\begin{array}{cc}p_{Y}& f\end{array}\right]}$ is a deflation it belongs to $\cj$-inj, and so
$\scriptsize{\left[\begin{array}{cc}p_{Y}\delta_{PY}& f\end{array}\right]}$ belongs to $\cj$-inj.

$\bullet$ \emph{Characterization of weak equivalences:} $f:X\ra Y$ is a morphism of $\cc$ such that
$\ul{\cc}(LN,\ol{f})$ is an isomorphism for all $N$, if and only $\ul{\cd}(N,\ol{Rf})$ is an isomorphism for all $N$, if and only if $\ul{\cc}(N,\cone(\ol{Rf}))=0$ for all $N$. For $N=\cone(\ol{Rf})$ we conclude that $\cone(\ol{Rf})=0$, \ie $Rf$ is a stable isomorphism.

$\bullet$ \emph{Characterizations of cofibrations:} use the factorization associated to $\ci$ and the retract argument to show that $\ci$-cof is included
in the class of inflations with cokernel of the form $LN$, for $N$ an object of $\cd$, and their retracts.

$\bullet$ \emph{Characterization of trivial cofibrations:} Let $g:A\ra B$ be an inflation with cokernel of the form $PLN$, and let $f$ be a morphism in $\cj$-inj. A commutative square of the form
\[\xymatrix{ A\ar[r]\ar[d]_{g} & C\ar[d]^f \\
B\ar[r] & D
}
\]
is isomorphic to one of the form
\[\xymatrix{ A\ar[r]\ar[d]_{\scriptsize{\left[\begin{array}{c}\id \\ 0\end{array}\right]}} & C\ar[d]^f \\
A\oplus PLN\ar[r] & D
}
\]
which admits an easy lifting map. On the other hand, by using the factorization associated to $\cj$ and the retract argument we see that $\cj$-cof is included
in the class of inflations with cokernel of the form $PLN$, for $N$ an object of $\cd$, and their retracts.

$\bullet$ \emph{Cofibrant objects:} It is clear that every object is fibrant. By Lemma \ref{Factorizations provided by adjunctions} we have that every object $LN$ is cofibrant. Now, if $M$ is a direct summand of $LN$ then $0\ra M$ is a retract of $0\ra LN$, and so it is a cofibration, \ie $M$ is cofibrant. Finally, if $M$ is cofibrant we have a lifting
\[\xymatrix{0\ar[r]\ar[d] & LRM\ar[d]^{\delta_{M}} \\
M\ar[r]_{\id}\ar@{.>}[ur] & M
}
\]
and so $M$ is a direct summand of $LRM$.

$\bullet$\emph{The t-structure:} Given any object $M$, the factorization of $M\ra 0$ associated to $\ci$ gives a conflation
\[M\arr{\beta}E\ra LRSM
\]
with $E\ra 0$ in $\ci$-inj, \ie  $E$ belongs to the full subcategory $\cn$ of $\ul{\cc}$ formed by the objects $M$ such that $\ul{\cc}(LN,M)=0$ for all $N$ in $\cd$. Then, there is a triangle in $\ul{\cc}$
\[M\ra E\ra LRSM\ra SM
\]
with $E$ in $\cn$ and $LRSM$ in the full subcategory $^{\bot}\cn$ of $\ul{\cc}$ formed by the objects $M'$ such that $\ul{\cc}(M',M)=0$ for all $M$ in $\cn$, and so $(^{\bot}\cn,\cn)$ is a t-structure on $\ul{\cc}$. Since the weak equivalences are the morphisms with cone in $\cn$, we have that $\op{Ho}\cc\simeq\!^{\bot}\cn$. It is clear that $\Tria(\{LN\}_{N\in\cd})\subseteq\!^{\bot}{\cn}$. On the other hand, if $M$ belongs to $\!^{\bot}\cn$, in the triangle we have
\[M\arr{0} E\ra LRSM\ra SM
\]
and so $LRSM\cong E\oplus SM$. But since $E\in\cn$, we deduce $E=0$ and so $SM\cong LRSM$, which implies that $M\in\Tria(\{LN\}_{N\in\cd})$.
\end{proof}

Notice that by applying the proposition above to $\cd^{op}$ we get the dual model structure in $\cd$, the dual t-structure and the dual description of the localization of $\cd$ with respect to the class $\cw_{\cd}$ of morphisms $g$ such that $L(g)$ is a stable isomorphism.

\begin{corollary}
Assume $R(\delta_{M})$ has injective kernel for each object $M$ of $\cc$ and $L(\eta_{N})$ has projective cokernel for each object $N$ of $\cd$.
\begin{enumerate}[1)]
\item The functors $L$ and $R$ induce mutually quasi-inverse triangle equivalences between $\cc[\cw_{\cc}^{-1}]$ and $\cd[\cw_{\cd}^{-1}]$.
\item If both $\cc$ and $\cd$ have finite direct and inverse limits, then they are model categories and the pair $(L,R)$ becomes a Quillen equivalence between $\cc$ and $\cd$.
\item If $\cc$ has kernels and $\cd$ has cokernels, then they are relevant derivable categories and the functors $L$ and $R$ induce mutually quasi-inverse equivalences between the corresponding derivators
\[\xymatrix{\mathbb{D}\cc\ar@<1ex>[d]^{\mathbb{R}R}_{\wr} \\
\mathbb{D}\cd\ar@<1ex>[u]^{\mathbb{L}L}
}
\]
\end{enumerate}
\end{corollary}
\begin{proof}
1) Let $\cn_{\cc}$ be the full subcategory of $\ul{\cc}$ formed by the objects $M$ such that $\ul{\cc}(LN,M)=0$ for all $N$ in $\cd$. Let $\cn_{\cd}$ be the full subcategory of $\ul{\cd}$ formed by the objects $N$ such that $\ul{\cd}(N,RM)=0$ for all $M$ in $\cc$. We know that $\cn_{\cc}$ is an aisle in $\ul{\cc}$ and that $\cc[\cw_{\cc}^{-1}]$ is the triangle quotient $\ul{\cc}/\cn_{\cc}$. Dually, we have that $\cn_{\cd}$ is a coaisle in $\ul{\cd}$ and that $\cd[\cw_{\cd}^{-1}]$ is the triangle quotient $\ul{\cd}/\cn_{\cd}$. Notice that $\cn_{\cc}$ is formed by all the objects $M$ such that $RM=0$ in $\ul{\cd}$, and that $\cn_{\cd}$ is formed by all the objects $N$ such that $LN=0$ in $\ul{\cc}$. Then, we have well defined adjoint triangle functors
\[\xymatrix{\ul{\cc}/\cn_{\cc}\ar@<1ex>[d]^{\ul{R}} \\
\ul{\cd}/\cn_{\cd}\ar@<1ex>[u]^{\ul{L}}
}
\]
Since the unit of the adjuntion $\eta_{N}:N\ra RLN$ becomes an isomorphism in $\ul{\cd}/\cn_{\cd}$ and the counit of the adjuntion $\delta_{M}:LRM\ra M$ becomes an isomorphism in $\ul{\cc}/\cn_{\cc}$, we have that these functors are mutually quasi-inverse.

2) Let us show that $(L,R)$ is a Quillen adjunction. To see that $L$ preserves cofibrations, thanks to the retract argument applied to
  the appropriate factorization for $\cd$, it suffices to show that $L\left(\scriptsize{\left[\begin{array}{c}\eta_{IX}\iota_{X}\\
          f\end{array}\right]}\right)$ is a cofibration. Since the cofibrations are closed under compositions it suffices to show that
  $L\left(\scriptsize{\left[\begin{array}{cc}\eta_{IX}&0\\  0&\id\end{array}\right]}\right)$ and
  $L\left(\scriptsize{\left[\begin{array}{c}\iota_{X}\\ f\end{array}\right]}\right)$ are cofibrations. Since
  $\scriptsize{\left[\begin{array}{c}\iota_{X}\\ f\end{array}\right]}$ is an inflation, then Lemma \ref{Factorizations provided by
    adjunctions} tells us that its image under $L$ is a cofibration. Now, since $L\left(\scriptsize{\left[\begin{array}{cc}\eta_{IX}&0\\
  0&\id\end{array}\right]}\right)$ is a section with projective cokernel direct summand of $LRLIX$, by Lemma \ref{Factorizations provided by adjunctions} it is cofibration. Dually, $R$ preserves fibrations.

  Let us show that $(L,R)$ is a Quillen equivalence. Let $f:LN\ra M$
  be a weak equivalence, where $N$ is an object of $\cd$ (and so
  cofibrant) and $M$ is an object of $\cc$ (and so fibrant). Since $R$
  is a right Quillen adjoint functor, it preserves weak equivalences
  between fibrant objects. Therefore, $\tau_{N,M}(f)=(Rf)\eta_{N}$ is
  a weak equivalence.
  
  3) We use Corollary \ref{A recognition criterion in Frobenius categories} to see that $\cc$ and $\cd$ are relevant derivable categories. Now, since $\cc$ and $\cd$ are derivable categories, then the associated prederivators are in fact derivators (\emph{cf.}\ Corollaire 2.28 of \cite{Cisinski2003}). It is easy to prove that $R$ is left exact in the sense of subsection 1.9 of \cite{Cisinski2003}. For instance, since $R$ has a left adjoint, then it preserves pullbacks. On the other hand, it is easy to prove that $R$ satisfies the ``approximation property'' of subsection 3.6 of \cite{Cisinski2003}. Therefore, by Theorem 3.12 of \cite{Cisinski2003} we know that the corresponding derived functor $\mathbb{R}R:\mathbb{D}\cc\ra\mathbb{D}\cd$ is an equivalence of derivators. Dually, $\mathbb{L}L:\mathbb{D}\cd\ra\mathbb{D}\cc$ is an equivalence of derivators.
\end{proof}

\begin{lemma}
\begin{enumerate}[a)]
\item If a morphism $f$ of $\cd$ is an inflation, then $L(f)\in\ci$-cof.
\item If $f$ is an inflation of $\cc$ with $\op{cok}(f)$ a direct summand of some LRM, then $f\in\ci$-cof.
\end{enumerate}
\end{lemma}
\begin{proof}
  $(a)$ We do the factorization associated to $\ci$ for $L(g)$. By the
  retract argument, it suffices to have that $L(g)$ has the left
  lifting property with respect to the composition
  $\scriptsize{\left[\begin{array}{cc}0&\id\end{array}\right]}\alpha$.
  But it is clear that $L(g)$ has the left lifting property with
  respect to
  $\scriptsize{\left[\begin{array}{cc}0&\id\end{array}\right]}$ and
  with respect to $\alpha$ (since $R(\alpha)$ is a retraction with
  injective kernel).

$(b)$ Consider the conflation $X\arr{f}Y\ra\op{cok}(f)$ with $\op{cok}(f)\oplus Q=LRM$. By adding the conflation
$0\ra Q\arr{\id}Q$ we get the conflation
\[\xymatrix{X\ar[r]^{\scriptsize{\left[\begin{array}{c}f\\ 0\end{array}\right]}\ \ \ } & Y\oplus Q\ar[r] & LRM.
}
\]
Hence, $f$ is a retract of and inflation with cokernel $LRM$.
\end{proof}

\begin{example}
Let $A$ be a unital dg algebra over $k$, and let $\cc A$ be the category of unital right dg $A$-modules.
As we will see in subsection \ref{The relative derived category}, Theorem \ref{Factorizations provided by adjunctions} provides a model structure on $\cc A$ whose associated homotopy category $\op{Ho}\cc A$ is the relative derived category $\cd_{rel}A$ (\emph{cf.}\ \cite{Keller1998a}).
\end{example}

\begin{example}
Let $\cc$ be a $k$-linear Frobenius category. By considering the adjuntion
\[\xymatrix{\cc\ar@<1ex>[d]^{\id} \\
\cc\ar@<1ex>[u]^{\id}
}
\]
we deduce from Theorem \ref{Factorizations provided by adjunctions} that $\cc$ admits two model structures according to which every object is fibrant and cofibrant, with the following classes of morphisms:
\begin{enumerate}
\item \emph{The projective model structure}:
\begin{enumerate}
\item The weak equivalences are the stable isomorphisms.
\item The cofibrations are the inflations and their retracts.
\item The fibrations are the morphisms with the right lifting property with respect to the morphisms of the form $0\ra P$ for some projective $P$.
\item The trivial cofibrations are the inflations with projective cokernel and their retracts.
\item The trivial fibrations are the morphisms with the right lifting property with respect to the morphisms of the form $S^{-1}M\ra PM$ for some $M$.
\end{enumerate}
\item \emph{The injective model structure}:
\begin{enumerate}
\item The weak equivalences are the stable isomorphisms.
\item The cofibrations are the morphisms with the left lifting property with respect to the morphisms of the form $I\ra 0$ for some injective $I$.
\item The fibrations are the deflations and their retracts.
\item The trivial cofibrations are the morphisms with the left lifting property with respect to the morphisms of the form $IM\ra SM$ for some $M$.
\item The trivial fibrations are the deflations with injective kernels and their retracts.
\end{enumerate}
\end{enumerate}
The localization with respect to these weak equivalences is triangle equivalent to the stable category $\ul{\cc}$. Hence, from the viewpoint of Homotopy Theory, things work nicely in Frobenius categories because they have model structures according to which every object is fibrant and cofibrant. See the preprint of D.-C. Cisinski \cite{Cisinski2003} for related considerations.
\end{example}
\bigskip

\section{$A_{[0,\infty[}$-algebras and their bar construction}\label{A zero infinity algebras and their bar construction}

An \emph{$A_{[0,\infty[}$-algebra over $k$}, also called \emph{curved $A_{\infty}$-algebra} or \emph{weak $A_{\infty}$-algebra}, is a graded $k$-module $A=\bigoplus_{p\in\Z}A^p$ together with a family of morphisms of graded $k$-modules (`multiplications')
\[m_{i}:A^{\otimes i}\ra A\ko i\geq 0\ko
\]
homogeneous of degree $|m_{i}|=2-i$, satisfying the identity
\[\sum_{j+k+l=p}(-1)^{jk+l}m_{i}(\id^{\otimes j}\otimes m_{k}\otimes\id^{\otimes l})=0\ko \text{for each } p\geq 0.
\]
Notice that $m_{0}$ is uniquely determined by the homogeneous element $m_{0}(1)$ of degree $2$ called the \emph{curvature} of $A$, and that $A_{\infty}$-algebras are precisely curved $A_{\infty}$-algebras with vanishing curvature. Following M. Kontsevich, we visualize the multiplication $m_{i}\ko i\geq 1$,  as a halfdisk whose upper arc is divided into segments, each of which symbolizes an `input', and whose base segment symbolizes the `output':
\[
\centeps{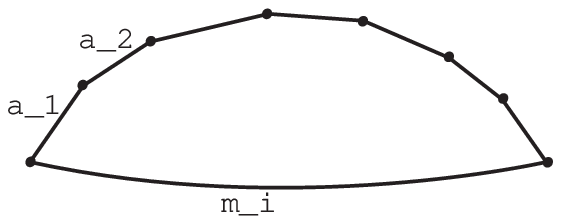}
\]
The morphism $\id^{\otimes j}\otimes m_{0}\otimes\id^{\otimes l}$ can be visualized as a bubble between positions $j$ and $j+1$:
\[
\centeps{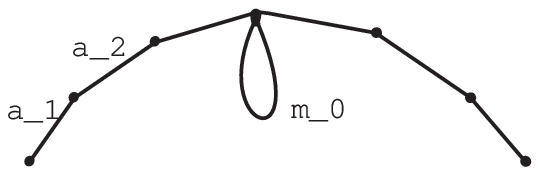}
\]
Using these representations, the defining identity is depicted as follows:
\[\sum \pm \begin{array}{c} \epsfbox{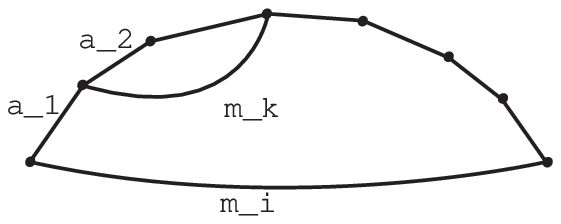}
\end{array} +
\]
\[+\sum\pm \begin{array}{c} \epsfbox{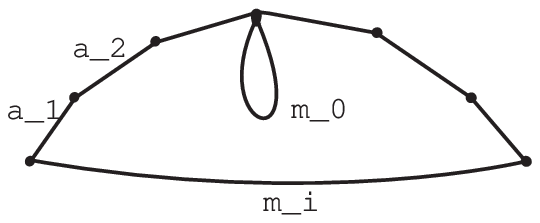}
\end{array}=0
\]
For instance, for $p=0$ we have $m_{1}m_{0}=0$, for $p=1$ we have $m_{1}m_{1}=m_{2}(m_{0}\otimes\id-\id\otimes m_{0})$ and for $p=2$ we have
$m_{3}(m_{0}\otimes\id^{\otimes 2}-\id\otimes m_{0}\otimes\id+\id^{\otimes 2}\otimes m_{0})+m_{1}m_{2}=m_{2}(m_{1}\otimes\id+\id\otimes m_{1})$. In particular, and in contrast with the situation when dealing with $A_{\infty}$-algebras, $m_{1}$ might not be a derivation and might not have vanishing square.
Thus, it is not clear how to define the `cohomology' of an $A_{[0,\infty[}$-algebra.

An $A_{[0,\infty[}$-algebra $A$ is \emph{strictly unital} if it is endowed with a homogeneous morphism $\eta: k\ra A$ of degree $0$, called the \emph{unit} of $A$, such that
\[m_{i}(\id_{A}\otimes\dots\otimes\id_{A}\otimes\eta\otimes\id_{A}\otimes\dots\otimes\id_{A})=0
\] 
for all $i\neq 2$ and
\[m_{2}(\id_{A}\otimes\eta)=m_{2}(\eta\otimes\id_{A})=\id_{A}.
\]

Let $(A,\{m_{i}\}_{i\geq 0})$ and $(A',\{m'_{i}\}_{i\geq 0})$ be two $A_{[0,\infty[}$-algebras over $k$. A \emph{morphism of $A_{[0,\infty[}$-algebras} is a sequence of morphisms of graded $k$-modules
\[f_{i}:A^{\otimes i}\ra A'\ko i\geq 1\ko
\]
homogeneous of degree $|f_{i}|=1-i$, satisfying the identity
\[\sum_{j+k+l=p}(-1)^{jk+l}f_{i}(\id^{\otimes j}\otimes m_{k}\otimes\id^{\otimes l})=\sum_{\underset{1\leq r\leq p}{i_{1}+\dots+i_{r}=p}}(-1)^{s}m'_{r}(f_{i_{1}}\otimes\dots\otimes f_{i_{r}})\ko
\]
for all $p\geq 0$, where
\[s=\sum_{2\leq u\leq r}\left((1-i_{u})\sum_{1\leq v\leq u-1}i_{v}\right).
\]
The right hand side of the equation is to be interpreted as $m'_{0}$ when $p=0$, and $s=1$ when $p=1$. Let $f=\{f_{i}\}_{i\geq 1}:A\ra A'$ and $f'=\{f'_{i}\}_{i\geq 1}:A'\ra A''$ be morphisms between $A_{[0,\infty[}$-algebras. The \emph{composition} is the morphism $f'f:A\ra A''$ with components
\[(f'f)_{p}=\sum_{\underset{1\leq r\leq p}{i_{1}+\dots+i_{r}=p}}(-1)^{s}f'_{r}(f_{i_{1}}\otimes\dots\otimes f_{i_{r}})\ko p\geq 1,
\]
where
\[s=\sum_{2\leq u\leq r}\left((1-i_{u})\sum_{1\leq v\leq u-1}i_{v}\right).
\]
A morphism $f$ of $A_{[0,\infty[}$-algebras is \emph{strict} if $f_{i}=0$ for $i\geq 2$.
For an $A_{[0,\infty[}$-algebra, the \emph{identity} of $A$ is the strict morphism $f:A\ra A$ given by $f_{1}=\id_{A}$. The category of $A_{[0,\infty[}$-algebras over $k$ will be denoted by $\op{Alg}_{[0,\infty[}$. Given an integer $n\geq 1$, the \emph{category of $A_{[0,n]}$-algebras over $k$}, denoted by
$\op{Alg}_{[0,n]}$, is the subcategory of $\op{Alg}_{[0,\infty[}$ formed by the algebras with multiplications $m_{i}=0$ for $i>n$ and morphisms $f_{i}=0$ for $i>n-1$.

\begin{example}
The objects of $\op{Alg}_{[0,2]}$ are the \emph{curved differential graded(=cdg) algebras} of \cite{Positselskii1993} (also called \emph{$Q$-algebras} sometimes). They are graded $k$-modules $A=\bigoplus_{p\in\Z}A^p$ together with a homogenous element $c$ of degree $2$, a morphism of graded $k$-modules $d:A\ra A$ homogeneous of degree $1$ called the \emph{predifferential}, and a morphism of graded $k$-modules $A\otimes A\ra A\ko (a,b)\mapsto ab$ homogeneous of degree $0$ called the \emph{multiplication}, satisfying
\begin{enumerate}[1)]
\item $d(c)=0$,
\item $d(d(a))=ca-ac$ for each $a\in A$,
\item $d(ab)=(da)b+(-1)^{|a|}a(db)$ for each $a$ (homogeneous) and $b$.
\end{enumerate}
The morphisms of $\op{Alg_{[0,2]}}$ are the morphisms of graded $k$-modules homogeneous of degree $0$ preserving the curvature, commuting with the predifferentials and preserving the multiplications. Thus, the morphisms of $\op{Alg}_{[0,2]}$ are instances of the morphisms between cdg algebras of \cite{Positselskii1993}. One can easily construct cdg algebras. Indeed, if $A$ is a graded algebra and $x$ is a homogeneous element of degree $1$, then we may set $m_{0}(1)= x^{2}\ko m_{1}(a)=xa-(-1)^{|a|}ax$ for $a$ homogeneous and $m_{2}$ equal to the product on $A$.
\end{example}

The bar construction of an $A_{\infty}$-algebra over $k$ yields a cocomplete augmented dg $k$-coalgebra \cite{Lefevre03}. In the presence of a non-zero curvature, the corresponding coalgebra is still a counital dg $k$-coalgebra, but the coaugmentation $\varepsilon$ is not compatible with the codifferential $d$ anymore, since one has $d\varepsilon\neq 0$. However, the coalgebra obtained is cocomplete augmented as a graded coalgebra. It means that it is an augmented graded coalgebra whose associated reduced coalgebra is cocomplete, \ie it is the colimit of the \emph{primitive filtration} formed by the kernels of the successive comultiplications. This leads to the category of \emph{cocomplete graded-augmented counital dg $k$-coalgebras}, $\op{CgaCdg}$, whose objects are $(C,d,\Delta,\varepsilon,\eta)$, such that $(C,d,\Delta,\eta)$ is a counital dg coalgebra and $(C,\Delta,\varepsilon, \eta)$ is a cocomplete augmented graded coalgebra, and the morphisms are morphisms of coalgebras compatible with the augmentations, the counits and the codifferentials.

Let $A$ be an $A_{[0,\infty[}$-algebra over $k$ and let $C$ be a cocomplete graded-augmented counital dg $k$-coalgebra. A straightforward calculation shows that we can endow the graded $k$-module $\HOM_{k}(C,A)$, whose $n$th component is formed by the $k$-linear morphisms $C\ra A$ homogeneous of degree $n$, with the curvature $b_{0}:=m_{0}^{A}\eta$, the first multiplication $b_{1}(f):=m_{1}^{A}f-(-1)^{|f|}fd_{C}$ and the following \emph{convolutions}
\[b_{n}(f_{1},\dots, f_{n}):=m_{n}^{A}(f_{1}\otimes\dots\otimes f_{n})\Delta^{(n)}\ko n\geq 2
\]
such that $\HOM_{k}(C,A)$ becomes an $A_{[0,\infty[}$-algebra over $k$. The set of \emph{twisting cochains} from $C$ to $A$ is the subset $\op{Tw}(C,A)$ of $\Hom^{1}_{k}(C,A)$ formed by the maps $\tau$ such that:
\begin{enumerate}
\item $\tau\varepsilon=0$,
\item $\tau$ satisfies the \emph{Maurer-Cartan equation} associated to $\HOM_{k}(C,A)$, \ie $\sum_{n\geq 0}b_{n}(\tau^{\otimes n})=0$.
\end{enumerate}

Notice that the sum above makes sense since $\tau$ is killed by the augmentation and $C$ is cocomplete graded-augmented.

\begin{proposition}
Given an $A_{[0,\infty[}$-algebra $A$ over $k$, the functor
\[\op{CgaCdg}\ra\op{Sets}\ko C\mapsto\op{Tw}(C,A)
\]
is representable. The \emph{bar construction} of $A$ is a representative $BA$. Moreover, the assignment $A\mapsto BA$ extends to a fully faithful covariant functor
\[B:\op{Alg_{[0,\infty[}}\ra\op{CgaCdg}
\]
such that the isomorphism $\op{Tw}(?,A)\cong\op{CgaCdg}(?,BA)$ is natural in $A$.
\end{proposition}
\begin{proof}
One checks that we can take $BA$ to be the tensor coalgebra $T^c(SA)=\bigoplus_{i\geq 0}(SA)^{\otimes i}$, such that its comultiplication `separates tensors'
\[\Delta(sx_{1},\dots, sx_{i})=\sum_{j=0}^{i}(sx_{1},\dots, sx_{j})\otimes(sx_{j+1},\dots, sx_{i}),
\]
where $(sx_{1},\dots, sx_{i})$ stands for $sx_{1}\otimes\dots\otimes sx_{i}$ and the empty parentheses $()$ are to be interpreted as $1_{k}$, endowed with a codifferential $d_{BA}$ which takes into account the multiplications $m_{i}\ko i\geq 0$ of $A$. More precisely, $d_{BA}$ is the unique coderivation of $BA$ such that the composition
\[BA\arr{d_{BA}}BA\arr{p_{0}}k
\]
vanishes (where $p_{0}$ is the projection on $k$), and the composition 
\[BA\arr{d_{BA}}BA\arr{p_{1}}SA
\]
(where $p_{1}$ is the projection on $SA$) has components $-sm_{i}w^{\otimes i}\ko i\geq 0$. It is straightforward to check that $m_{i}\ko i\geq 0$, define a structure of $A_{[0,\infty[}$-algebra on $A$ if and only if $d^2_{BA}=0$. The bijection $\op{CgaCdg}(C,BA)\arr{\sim}\op{Tw}(C,A)$ takes $F$ to the map $\tau$ given by the composition
\[C\arr{F}BA\arr{p_{1}}SA\arr{w}A.
\]
Observe that $\tau\varepsilon=0$. It is straightforward to check that $F$ is compatible with the codifferentials if and only if $\tau$ is a solution of the Maurer-Cartan equation of the $A_{[0,\infty[}$-algebra $\HOM_{k}(C,A)$. To define $B$ on morphisms we do the following. Let $A$ and $A'$ be two $A_{[0,\infty[}$-algebras, and let $f_{i}:A^{\otimes i}\ra A'\ko i\geq 1$, be a family of morphisms of graded $k$-modules of degree $|f_{i}|=1-i$. Consider the morphism of coaugmented graded coalgebras $Bf: BA\ra BA'$ which is compatible with the counits and the coaugmentations and takes into account the morphisms $f_{i}\ko i\geq 1$. More precisely, $Bf$ is the unique morphism of graded coalgebras which is compatible with the counits and the coaugmentations and such that its associated morphism of graded coalgebras $\ol{Bf}:\ol{BA}\ra\ol{BA'}$ between the corresponding reduced tensor coalgebras satisfies
\[\ol{p_{n}}\ol{Bf}=\scriptsize{\left[\begin{array}{ccc}F_{1}&F_{2}&\dots\end{array}\right]}^{\otimes n}\ol{\Delta}^{(n)}\ko n\geq 1,
\]
where $\ol{p_{n}}:\ol{BA}\ra(SA)^{\otimes n}$ is the projection, $\ol{\Delta}$ is the `separating tensors' comultiplication of $\ol{BA}$ and $F_{i}=sf_{i}w^{\otimes i}\ko i\geq 1$. It is straightforward to check that $f_{i}\ko i\geq 1$, define a morphism of $A_{[0,\infty[}$-algebras if and only if $Bf$ commutes with the codifferentials of $BA$ and $BA'$.
\end{proof}

It is clear from the proof of the proposition that the $A_{[0,\infty[}$-algebra structures on a graded $k$-module $A$ are in bijection with the coderivations $d$ of the graded coalgebra $(T^{c}(SA),\Delta)$ making $(T^c(SA),d,\Delta, \varepsilon, \eta)$ into a cocomplete graded-augmented counital dg $k$-coalgebra.
\bigskip

\section{Bar/cobar adjunction for cdg algebras}\label{Bar cobar adjunction for cdg algebras}

\begin{proposition}
Given a cocomplete graded-augmented counital dg $k$-coalgebra $C$, the functor
\[\op{Alg}_{[0,2]}\ra\op{Sets}\ko A\mapsto\op{Tw}(C,A)
\]
is corepresentable. The \emph{cobar construction} of $C$ is a corepresentative $\Omega C$. Moreover, the assignment $C\mapsto\Omega C$ extends to a covariant functor
\[\Omega:\op{CgaCdg}\ra\op{Alg}_{[0,2]}
\]
such that the isomorphism $\op{Tw}(C,?)\cong\op{Alg}_{[0,2]}(\Omega C,?)$ is natural in $C$.
\end{proposition}
\begin{proof} 
Let $\pi: C \to \ol{C}$ be the cokernel of $\eps$ and $\rho : \ol{C} \to C$ its canonical section whose image is the kernel of the counit $\eta$. Notice that $d$ does not induce a map in $\ol{C}$. We put $\ol{d}=\pi d\rho$. We let $\ol{\Delta}$ be the associative comultiplication of $\ol{C}$ induced by $\Delta$.  One
checks that one can take $\Omega C$ to be the reduced tensor algebra $\ol{T}(S^{-1}\ol{C})=\bigoplus_{i\geq 1}(S^{-1}\ol{C})^{\otimes i}$, endowed with a differential which takes into account the maps $\ol{d}:\ol{C}\ra\ol{C}$ and $\ol{\Delta}:\ol{C}\ra\ol{C}\otimes\ol{C}$ induced by the codifferential $d$ and the comultiplication $\Delta$ of $C$. The curvature of $\Omega C$ is given by the composition
\[wd\varepsilon: k\ra S^{-1}\ol{C}\ko
\]
where $w$ is the degree shift morphism.
\end{proof}

As a consequence, we get
\begin{corollary}
The bar and cobar constructions
\[\xymatrix{\op{Alg}_{[0,2]}\ar@<1ex>[d]^{B} \\
\op{CgaCdg}\ar@<1ex>[u]^{\Omega}
}
\]
form a pair of adjoint functors.
\end{corollary}
\bigskip

Given an $A_{[0,\infty[}$-algebra $A$ we have the following isomorphisms
\[\op{Alg}_{[0,2]}(\Omega BA,\Omega BA)\cong\op{CgaCdg}(BA, B\Omega BA)\underset{\sim}{\stackrel{B}{\leftarrow}}\op{Alg}_{[0,\infty[}(A,\Omega BA)
\]
Therefore, the identity morphism $\id_{\Omega BA}$ corresponds to a morphism of $A_{[0,\infty[}$-algebras from $A$ to the so-called \emph{envelopping cdg algebra} of $A$
\[A\ra \Omega BA.
\]
Moreover, this morphism is universal among the morphisms of $A_{[0,\infty[}$-algebras from $A$ to a curved dg algebra. We expect this to be a homotopy equivalence of precomplexes (\ie `complexes' such that the square of the differential is not necessarily zero) over $k$, at least when $A$ is an augmented $A_{[0,\infty[}$-algebra.
\bigskip

\section{$A_{[0,\infty[}$-modules and their bar construction}\label{A zero infinity modules and their bar construction}

\subsection{Basic notions}\label{Basic notions}
Let $A$ be an $A_{[0,\infty[}$-algebra over $k$. A \emph{right $A_{[0,\infty[}$-module} over $A$ is a graded $k$-module $M=\bigoplus_{p\in\Z}M^p$ endowed with a family of morphisms of graded $k$-modules
\[m_{i}^M:M\otimes A^{\otimes(i-1)}\ra M\ko i\geq 1,
\]
homogeneous of degree $|m^M_{i}|=2-i$ satisfying the identity
\[\sum_{j+k+l=p}(-1)^{jk+l}m_{i}^M(\id^{\otimes j}\otimes m_{k}\otimes\id^{\otimes l})=0\ko p\geq 1.
\]
For $j>0$ the term $m^M_{i}(\id^{\otimes j}\otimes m_{k}\otimes\id^{\otimes l})$ should be interpreted as
\[m^M_{i}(\id^{\otimes j}\otimes m_{k}^{A}\otimes\id^{\otimes l}):M\otimes A^{\otimes (p-1)}\ra M,
\]
and for $j=0$ as
\[m_{i}^M(m_{k}^M\otimes\id^{\otimes l}):M\otimes A^{\otimes(p-1)}\ra M.
\]
For instance, for $p=1$ we have $m_{1}^Mm_{1}^M=-m_{2}^M(\id\otimes m_{0}^M)$, and for $p=2$ we have $m_{1}^Mm_{2}^M=m_{2}^M(m_{1}^M\otimes\id+\id\otimes m_{1}^{A})+m_{3}^M(\id\otimes m_{0}^{A}\otimes\id+\id\otimes\id\otimes m_{0}^{A})$. In particular, $A$ is not an $A_{[0,\infty[}$-module over itself with the regular structure, $m_{1}^M$ is not a derivation and it is not clear what the cohomology of an $A_{[0,\infty[}$-module should be.

If $A$ is strictly unital with unit $\eta$, an $A_{[0,\infty[}$-module $M$ over $A$ is \emph{strictly unital} if $m^M_{i}(\id_{M}\otimes\dots\otimes\id\otimes\eta\otimes\id\otimes\dots\otimes\id)=0$ for all $i\geq 3$ and
\[m^{M}_{2}(\id_{M}\otimes\eta)=\id_{M}.
\]

Let $M$ and $N$ be two $A_{[0,\infty[}$-modules over $A$. A \emph{morphism of $A_{[0,\infty[}$-modules} from $M$ to $N$ is a family of morphisms of graded $k$-modules
\[f_{i}:M\otimes A^{\otimes(i-1)}\ra N\ko i\geq 1,
\]
homogeneous of degree $|f_{i}|=1-i$ satisfying the identity
\[\sum_{j+k+l=p}(-1)^{jk+l}f_{i}(\id^{\otimes j}\otimes m_{k}\otimes\id^{\otimes l})=\sum_{r+s=p}m_{s+1}^N(f_{r}\otimes\id^{\otimes s})\ko p\geq 1.
\]
Notice that we can not have $j=k=0$. A morphism $f$ is \emph{strict} if $f_{i}=0$ for all $i\neq 1$. A morphism $f:M\ra N$ between strictly unital $A_{[\infty,0[}$-modules is \emph{strictly unital} if
\[f_{i}(\id_{M}\otimes\id\otimes\dots\otimes\id\otimes\eta\otimes\id\otimes\dots\otimes\id)=0\ko i\geq 2.
\]
If $M\ko N$ and $T$ are three $A_{[0,\infty[}$-modules over $A$, and $g:M\ra N$ and $f:N\ra T$ are two morphisms of $A_{[0,\infty[}$-modules, the \emph{composition} $fg$ is defined by the family
\[(fg)_{p}=\sum_{k+l=p}f_{1+l}(g_{k}\otimes\id^{\otimes l})\ko p\geq 1.
\]

Let $f\ko g:M\ra N$ be two morphisms of $A_{[0,\infty[}$-modules over $A$. A \emph{homotopy of $A_{[0,\infty[}$-modules from $f$ to $g$} is a family of morphisms of graded $k$-modules
\[h_{i}:M\otimes A^{\otimes(i-1)}\ra N\ko i\geq 1,
\]
homogeneous of degree $|h_{i}|=-i$ satisfying the identity
\[f_{p}-g_{p}=\sum_{r+s=p}(-1)^{s}m_{1+s}^{N}(h_{r}\otimes\id^{\otimes s})+\sum_{j+k+l=p}(-1)^{jk+l}h_{i}(\id^{\otimes j}\otimes m_{k}\otimes\id^{\otimes l})\ko p\geq 1
\]
If $f$ and $g$ are two strictly unital morphisms of $A_{[\infty,0[}$-modules, a homotopy $h$ between $f$ and $g$ is \emph{strictly unital} if
\[h_{i}(\id_{M}\otimes\id\otimes\dots\otimes\eta\otimes\id\otimes\dots\otimes\id)=0\ko i\geq 2.
\]
Two morphisms of $A_{[0,\infty[}$-modules $f$ and $g$ are \emph{homotopic} if there exists a homotopy between $f$ and $g$.
A morphism $f$ of $A_{[0,\infty[}$-modules is \emph{null-homotopic} if $f$ and $0$ are homotopic.
An $A_{[0,\infty[}$-module $M$ is \emph{contractible} if $\id_{M}$ is null-homotopic.

We denote by $\op{Nod}_{[0,\infty[}A$ the category of all the right $A_{[0,\infty[}$-modules over $A$, and by $\op{Nod}^{strict}_{[0,\infty[}A$ the subcategory formed by all the right $A_{[0,\infty[}$-modules over $A$ with strict morphisms. (The `$\op{N}$' in $\op{Nod}A$ comes from `\textbf{n}on-unital', and goes back to Kenji Lef{\`e}vre-Hasegawa's thesis \cite{Lefevre03}.) Also, if $A$ is a strictly unital $A_{\infty}$-algebra, then $\Mod_{[\infty,0[}A$ is the category of strictly unital $A_{[\infty,0[}$-modules over $A$ with strictly unital morphisms.

If $(M,m_{1}^M,m_{2}^M,\dots)$ is an $A_{[0,\infty[}$-module over $A$, then the shifted graded $k$-module $SM$ inherits naturally a structure of right $A_{[0,\infty[}$-module over $A$ with multiplications $m_{i}^{SM}:=(-1)^{i}m_{i}^{M}(w\otimes\id^{\otimes(i-1)})$. Moreover, the shift easily extends to an automorphism
\[S:\op{Nod}^{strict}_{[0,\infty[}A\arr{\sim}\op{Nod}^{strict}_{[0,\infty[}A.
\]

Given an $A_{[0,\infty[}$-algebra $A$, we denote by $\op{MC}(A)$ the set of elements $a\in A^1$ such that $m_{i}^{A}(a^{\otimes i})=0$ for $i\gg 0$ and $a$ satisfies the \emph{Maurer-Cartan equation} associated to $A$:
\[\sum_{i\geq 0}m_{i}^{A}(a^{\otimes i})=0.
\] 
Since the world of modules is a `linearization' of the world of algebras, we do not have the concept of the Maurer-Cartan equation of an $A_{[0,\infty[}$-module, but rather the concept of the tangent space to such a (phantom) equation. Indeed, let $M$ be a right $A_{[0,\infty[}$-module over $A$ and let $a$ be an element of $\op{MC}(A)$ such that for each $m\in M$ there exists $i_{0}\geq 1$ satisfying $m_{i}^M(m\otimes a^{\otimes(i-1)})=0$ for all $i\geq i_{0}$.
We define the \emph{tangent space in $a$ to the Maurer-Cartan equation of $M$} to be the $k$-submodule $T_{a}\op{MC}(M)$ of $M^1$ formed by the elements $m\in M^1$ such that
\[\sum_{i\geq 1}m_{i}^{M}(m\otimes a^{\otimes(i-1)})=0.
\]

\subsection{Twisting cochains}\label{Twisting cochains}
Now we need a technical lemma which will be very useful for understanding the tangent space to the Maurer-Cartan equation of an $A_{[0,\infty[}$-module (\emph{cf.}\ Proposition \ref{Twisting cochains}) and so, for the bar construction of an $A_{[0,\infty[}$-module.

\begin{lemma}
Let $C$ be a cocomplete graded-augmented counital dg $k$-coalgebra, $A$ an $A_{[0,\infty[}$-algebra, $\tau$ a twisting cochain from $C$ to $A$ and $M$ an $A_{[0,\infty[}$-module over $A$. Then the map $g:M\otimes C\ra M$ defined by
\[g:=\sum_{i\geq 0}(-1)^{i}m_{i+1}^{M}(\id_{M}\otimes \tau^{\otimes i}\Delta^{(i)})
\]
(with $\Delta^{(0)}:=\eta$ and $\Delta^{(1)}:=\id$) satisfies the identity
\[g(\id_{M}\otimes d_{C})+g(g\otimes\id_{C})(\id_{M}\otimes\Delta)=0.
\]
\end{lemma}
\begin{proof}
By taking into account that $\tau$ is killed by the augmentation and that $C$ is a cocomplete graded-augmented coalgebra, \ie that $C=k\oplus(\underset{n\geq 2}{\op{colim}}\ker\ol{\Delta}^{(n)})$, one can check that for each element $m\otimes c$ of $M\otimes C$ the expression $g(m\otimes c)$ makes sense since it becomes a finite sum. It suffices to prove that the expression above vanishes when applied to an element of the form $m\otimes 1$ or $m\otimes x$ with $x\in\ker\ol{\Delta}^{(n+1)}$ for some $n\geq 1$. The first case is trivial, and in the second case, one can prove that
\[(g(\id_{M}\otimes d_{C})+g(g\otimes\id_{C})(\id_{M}\otimes\Delta))(m\otimes x)
\]
equals the following sum of $n$ sums
\[\sum_{t=1}^{n}\left((-1)^{s(t)}\sum_{j+k+l=t+1}(-1)^{jk+l}m_{i}^M(\id^{\otimes j}\otimes m_{k}\otimes\id^{\otimes l})(m\otimes\tau x_{t1}\otimes\dots\otimes\tau x_{tt})\right)=0,
\]
where we assume $\ol{\Delta}^{(t)}(x)=x_{t1}\otimes\dots\otimes x_{tt}$ for the sake of simplicity, and where
\[s(t)=t(m+1)+\sum_{2\leq u\leq t}\sum_{1\leq v\leq u-1}|x_{tv}|.
\]
\end{proof}

The next proposition is the $A_{[0,\infty[}$-counterpart of \cite[Lemma 3.4.1]{CKTB}, which allows us to describe the tangent space in $a$ to the Maurer-Cartan equation of an $A_{[0,\infty[}$-module as the space of the $1$-cocycles of a certain complex defined by twisting with $a$.

\begin{proposition}
Let $A$ be an $A_{[0,\infty[}$-algebra, $M$ an $A_{[0,\infty[}$-module over $A$ and $a$ an element of $\op{MC}(A)$ such that for each $m\in M$ there exists $i_{0}\geq 1$ satisfying $m_{i}^M(m\otimes a^{\otimes(i-1)})=0$ for all $i\geq i_{0}$.
\begin{enumerate}[1)]
\item The graded $k$-module $M$ becomes a complex over $k$ endowed with the map
\[d_{a}(m):=\sum_{i\geq 1}(-1)^{(i-1)(|m|+1)}m_{i}^M(m\otimes a^{\otimes(i-1)}).
\]
Moreover, $(M,d_{a}^M)$ underlies a functor
\[F_{a}:\op{Nod}^{strict}_{[0,\infty[}A\ra \cc k
\]
compatible with the shift.
\item $m\in M^1$ belongs to $T_{a}\op{MC}(M)$ if and only if $d_{a}(m)=0$.
\item For $m\in M^0$, the map $a\mapsto d_{a}(m)$ yields a `vector field' on $\op{MC}(A)$, \ie $d_{a}(m)\in T_{a}\op{MC}(M)$ for all $a\in\op{MC}(A)$.
\end{enumerate}
\end{proposition}
\begin{proof}
If we apply Lemma \ref{Twisting cochains} to the coalgebra $C=BA$ we have 
\[d_{a}(m)=\sum_{i\geq 0}g(m\otimes (sa)^{\otimes i}),
\] 
and so
\begin{align}
d_{a}(d_{a}(m))=\sum_{p\geq 0}g(g\otimes\id_{BA})(\id_{M}\otimes\Delta)(m\otimes(sa)^{\otimes p})= \nonumber \\
=-\sum_{p\geq 0}g(\id_{M}\otimes d_{BA})(m\otimes(sa)^{\otimes p})=(-1)^{|m|+1}g\left(m\otimes\sum_{p\geq 0}d_{BA}((sa)^{\otimes p})\right). \nonumber
\end{align}
But, taking into account that $|sa|=1+|a|=2$, we have
\begin{align}
d_{BA}((sa)^{\otimes p})=-\sum_{j+k+l=p}(sa)^{\otimes j}\otimes(sm_{k}^{A}w^{\otimes k})((sa)^{\otimes k})\otimes(sa)^{\otimes l}= \nonumber \\
=-\sum_{j+k+l=p}(sa)^{\otimes j}\otimes sm_{k}^{A}(a^{\otimes k})\otimes(sa)^{\otimes l}, \nonumber
\end{align}
and so
\begin{align}
d_{a}(d_{a}(m))=(-1)^{|m|}\sum_{p\geq 0}\sum_{j+k+l=p}g(m\otimes(sa)^{\otimes j}\otimes sm_{k}^{A}(a)\otimes (sa)^{\otimes l})= \nonumber \\
=(-1)^{|m|}\sum_{j,l\geq 0}g(m\otimes(sa)^{\otimes j}\otimes\left(s\sum_{k\geq 0}m_{k}^{A}(a)\right)\otimes (sa)^{\otimes l})=0, \nonumber
\end{align}
since $a\in\op{MC}(A)$. The rest of the proof is straightforward.
\end{proof}

\subsection{The bar construction}\label{The bar construction}

Given an $A_{[0,\infty[}$-algebra $A$ over $k$, an $A_{[0,\infty[}$-module $M$ over $A$, a cocomplete graded-augmented counital dg $k$-coalgebra $C$ and a dg right $C$-comodule $N$, we can endow the graded $k$-module $\HOM_{k}(N,M)$ with the multiplications
\[\mu_{1}(f):=m_{1}^Mf-(-1)^{|f|}fd_{N}
\]
and
\[\mu_{i}(f\otimes\alpha_{1}\otimes\dots\otimes\alpha_{i-1}):=m_{i}^M(f\otimes\alpha_{1}\otimes\dots\otimes\alpha_{i-1})\Delta_{N}^{(i)}\ko i\geq 2,
\]
so that $\HOM_{k}(N,M)$ becomes a right $A_{[0,\infty[}$-module over $\HOM_{k}(C,A)$ regarded as an $A_{[0,\infty[}$-algebra with the multiplications $b_{i}\ko i\geq 0$, defined in section \ref{A zero infinity algebras and their bar construction}.

Let $\tau$ be a twisting cochain from $C$ to $A$. Notice that in general it is not true that $b_{i}(\tau^{\otimes i})=0$ for $i\gg 0$ and that for each $f$ in $\HOM_{k}(N,M)$ there exists a natural number $i_{0}\geq 1$ such that $\mu_{i}(f\otimes\tau^{(i-1)})=0$ for $i\geq i_{0}$. However, for each $c\in C$ we do have $b_{i}(\tau^{\otimes})(c)=0$ for $i\gg 0$, and for each $f\in\HOM_{k}(N,M)$ and $n\in N$ we do have $\mu_{i}(f\otimes\tau^{(i-1)})(n)=0$ for $i\gg 0$. Therefore, in practice it still makes sense to speak of the tangent space in $\tau$ to the Maurer-Cartan equation of the $A_{[0,\infty[}$-module $S^{-1}\HOM_{k}(N,M)$, denoted by $T_{\tau}\op{MC}(S^{-1}\HOM_{k}(N,M))$. By using the techniques of Proposition \ref{Twisting cochains} we can check that this tangent space identifies with the space of $0$-cocycles of the cochain complex $F_{\tau}(\HOM_{k}(N,M))$, consisting of the graded module $\HOM_{k}(N,M)$  together with the differential given by
\[d_{\tau}(f):=\sum_{i\geq 1}(-1)^{(i-1)(|f|+1)}\mu_{i}(f\otimes\tau^{\otimes(i-1)}).
\]

To find the bar construction over $C$ of $M$, we need to consider the category of counital dg right $C$-comodules, $\op{Com}C$. It is a Frobenius category since it is the category of $0$-cocycles of a certain exact dg category $\op{Com}_{dg}C$, \emph{cf.}\ \cite{Keller1999, Keller2006b}. The conflations are the short exact sequences which split in the category of graded $C$-comodules, and an object $N$ of $\op{Com}C$ is injective-projective if and only if its identity morphism $\id_{N}$ is \emph{null-homotopic}, \ie $\id_{N}=d_{N}h+hd_{N}$ for some morphism $h$ of graded $C$-comodules homogeneous of degree $-1$.
The corresponding stable category, denoted by $\ch C$ and called the \emph{category of counital $C$-comodules up to homotopy}, is the quotient of $\op{Com} C$ by the ideal of the null-homotopic morphisms. We also need the category of counital graded right $C$-comodules, $\op{Grmod}C$, with homogeneous morphisms of degree $0$.

Notice that the counital graded right $C$-comodule, $(M\otimes_{k}C, \id_{M}\otimes\Delta_{C})$, becomes a counital dg right $C$-comodule, $M\otimes_{\tau}C$, with the codifferential induced \cite[Lemme 2.1.2.1]{Lefevre03} by the map of the Lemma \ref{Twisting cochains}. Now, the isomorphisms
of $k$-modules
\[(\op{Grmod}C)(N[-p],M\otimes_{k}C)\arr{\sim}\Hom^p_{k}(N,M)\ko f\mapsto (\id_{M}\otimes\eta)f,
\]
induce an isomorphism of complexes over $k$
\[(\op{Com}_{dg}C)(N,M\otimes_{\tau}C)\arr{\sim}F_{\tau}(\HOM_{k}(N,M)),
\]
and the isomorphism between the $0$-cocycles gives us
\[(\op{Com}C)(N,M\otimes_{\tau}C)\arr{\sim}T_{\tau}\op{MC}(S^{-1}\HOM_{k}(N,M)).
\]
Therefore, we have

\begin{proposition}
Let $A$ be an $A_{[0,\infty[}$-algebra, $M$ an $A_{[0,\infty[}$-module over $A$, $C$ a cocomplete graded-augmented counital dg coalgebra and $\tau$ a twisting cochain from $C$ to $A$. The functor
\[\op{Com}C\ra\Mod k\ko N\mapsto T_{\tau}\op{MC}(S^{-1}\HOM_{k}(N,M))
\]
is representable. The \emph{bar construction} over $C$ of $M$ is a representative, denoted by $M\otimes_{\tau}C$. The assignment $M\mapsto M\otimes_{\tau}C$ extends to a functor
\[?\otimes_{\tau}C:\op{Nod}_{[0,\infty[}A\ra\op{Com}C\ko M\mapsto M\otimes_{\tau}C
\]
such that the isomorphism $(\op{Com}C)(?,M\otimes_{\tau}C)\arr{\sim}T_{\tau}\op{MC}(S^{-1}\HOM_{k}(?,M))$ is natural in $M$.
\end{proposition}

Given an $A_{[0,\infty[}$-algebra $A$, the composition of the projection with the shift to the right, $\tau: BA\arr{p_{1}}SA\arr{w}A$, is a twisting cochain from $BA$ to $A$. By using the Lemma \ref{Twisting cochains} one can prove that the structures of $A_{[0,\infty[}$-module over $A$ of a graded $k$-module $M$ are in bijection with the codifferentials making the counital graded comodule $M\otimes_{k}BA$ into a counital dg comodule. Similarly, the morphisms of $A_{[0,\infty[}$-modules from $M$ to $N$ are in bijection with the morphisms of dg $BA$-comodules from $M\otimes_{\tau}BA$ to $N\otimes_{\tau}BA$, and the homotopies of $A_{[0,\infty[}$-modules from $f$ to $g$ are in bijection with the homotopies of dg $BA$-comodules from the morphism induced by $f$ to the morphism induced by $g$.
\bigskip

\section{Bar/cobar adjunction for cdg modules}\label{Bar cobar adjunction for cdg modules}

Given an $A_{[0,n]}$-algebra, $A$, the category of \emph{$A_{[0,n]}$-modules over $A$} is the subcategory $\op{Nod}_{[0,n]}A$ of $\op{Nod}_{[0,\infty[}$ formed by the $A_{[0,\infty[}$-modules $M$ with multiplications $m_{i}^M=0$ for $i>n$ and morphisms $f_{i}=0$ for $i>n-1$. For $n=2$ we get the so-called \emph{curved dg(=cdg) modules} over a cdg algebra.

Let us consider a fixed (strictly) unital cdg algebra $A$ and denote by $\cc A$ the category of unital cdg $A$-modules. $\cc A$ is the category of $0$-cocycles of the exact dg category (\emph{cf.}\ \cite{Keller1999, Keller2006b}) $\cc_{dg}A$ whose objects are the unital cdg $A$-modules and whose morphisms are given by complexes of $k$-modules $\cc_{dg}(A)(L,M)$ with $n$th component formed by the morphisms of graded $k$-modules homogeneous of degree $n$ and with differential given by the commutator $d(f)=d_{M}f-(-1)^{|f|}fd_{L}$. Thus, as in the case of a dg algebra, $\cc A$ has a structure of a $k$-linear Frobenius category whose stable category is the \emph{category of unital cdg $A$-modules up to homotopy}, $\ch A$, defined precisely as in the case of unital dg algebras (\emph{cf.}\ the first example of subsecction \ref{Small object argument in Frobenius categories}).

Given a counital dg coalgebra $C$ we still have a cdg algebra $\HOM_{k}(C,A)$ as before (\emph{cf.}\ section \ref{A zero infinity algebras and their bar construction}). A solution $\tau$ of the Maurer-Cartan equation of this cdg algebra will suffice to `twist the cochain' since almost all the multiplications of $A$ are zero, and so we do not need any extra assumption on $\tau$.

\subsection{Bar construction for unital cdg modules}\label{Bar construction for unital cdg modules}
For a unital cdg $A$-module $M$, the counital graded $C$-comodule $(M\otimes_{k}C,\id_{M}\otimes\Delta_{C})$ becomes a counital dg $C$-comodule, $M\otimes_{\tau}C$, with the codifferential
\[d_{M\otimes_{\tau}C}:=d_{M}\otimes\id_{C}+\id_{M}\otimes d_{C}-(m_{2}^M\otimes\id)(\id\otimes\tau\otimes\id)(\id\otimes\Delta).
\]
As in subsection \ref{The bar construction}, the isomorphism of complexes over $k$
\[(\op{Com}_{dg}C)(N,M\otimes_{\tau}C)\arr{\sim}F_{\tau}(\HOM_{k}(N,M))
\]
induces an isomorphism between the $0$-cocycles which gives us
\[(\op{Com}C)(N,M\otimes_{\tau}C)\arr{\sim}T_{\tau}\op{MC}(S^{-1}\HOM_{k}(N,M)).
\]

\begin{proposition}
Let $A$ be a cdg algebra, $M$ a unital cdg $A$-module, $C$ a counital dg coalgebra and $\tau$ a solution of the Maurer-Cartan equation of the cdg algebra $\HOM_{k}(C,A)$. The functor
\[\op{Com}C\ra\Mod k\ko N\mapsto T_{\tau}\op{MC}(S^{-1}\HOM_{k}(N,M))
\]
is representable. The \emph{bar construction} over $C$ of $M$ is a representative, denoted by $M\otimes_{\tau}C$. Moreover, it can be chosen to yield a functor
\[?\otimes_{\tau}C:\cc A\ra\op{Com}C\ko M\mapsto M\otimes_{\tau}C
\]
such that the isomorphism $(\op{Com}C)(?,M\otimes_{\tau}C)\arr{\sim}T_{\tau}\op{MC}(S^{-1}\HOM_{k}(?,M))$ is natural in $M$.
\end{proposition}

\subsection{Cobar construction for counital dg comodules}\label{Cobar construction for counital dg comodules}

For a counital dg $C$-module $N$, the unital graded $A$-module $(N\otimes_{k}A,\id_{M}\otimes m_{2}^{A})$ becomes a unital cdg $A$-module, $N\otimes_{\tau}A$, with the differential
\[d_{N\otimes_{\tau}A}:=d_{N}\otimes\id_{A}+\id_{N}\otimes m_{1}^{A}+(\id\otimes m_{2}^A)(\id\otimes\tau\otimes\id)(\Delta_{N}\otimes\id).
\]
Now, the isomorphisms of $k$-modules
\[(\op{Grmod}A)(N\otimes_{k}A,M[p])\arr{\sim}\Hom^p_{k}(N,M)\ko f\mapsto f(\id_{N}\otimes\eta_{A})
\]
where $\op{Grmod} A$ is the category of graded unital $A$-modules and $\eta_{A}$ is the unit of $A$, induce an isomorphism of complexes over $k$
\[(\cc_{dg}A)(N\otimes_{\tau}A,M)\arr{\sim}F_{\tau}(\HOM_{k}(N,M))
\]
which yields an isomorphism between the $0$-cocycles
\[(\cc A)(N\otimes_{\tau}A,M)\arr{\sim}T_{\tau}\op{MC}(S^{-1}\HOM_{k}(N,M)).
\]

\begin{proposition}
Let $A$ be a cdg algebra, $C$ a counital dg coalgebra, $N$ a counital dg $C$-comodule and $\tau$ a solution of the Maurer-Cartan equation of the cdg algebra $\HOM_{k}(C,A)$. The functor
\[\cc A\ra\Mod k\ko M\mapsto T_{\tau}\op{MC}(S^{-1}\HOM_{k}(N,M))
\]
is corepresentable. The \emph{cobar construction} over $A$ of $N$ is a corepresentative, denoted by $N\otimes_{\tau}A$. Moreover, it can be chosen to yield a functor
\[?\otimes_{\tau}A:\op{Com}C\ra\cc A\ko N\mapsto N\otimes_{\tau}A
\]
such that the isomorphism $(\cc A)(N\otimes_{\tau}A,?)\arr{\sim}T_{\tau}\op{MC}(S^{-1}\HOM_{k}(N,?))$ is natural in $N$.
\end{proposition}

\subsection{The adjunction between Frobenius categories}\label{The adjunction between Frobenius categories}

Let $A$ be a unital cdg algebra, and $C$ a counital dg coalgebra. Let us consider the category $\cc A$ of unital cdg right $A$-modules, which is a Frobenius category as explained at the beginning of section \ref{Bar cobar adjunction for cdg modules}. Also, let us consider the category $\op{Com}C$ of counital dg right $C$-comodules, which is a Frobenius category as explained in subsection \ref{The bar construction}.

\begin{corollary}
Let $\tau$ be a solution of the Maurer-Cartan equation in $\HOM_{k}(C,A)$. Then the bar and cobar constructions yield a pair of adjoint exact functors
\[\xymatrix{ \cc A\ar@<1ex>[d]^{?\otimes_{\tau}C=:R} \\
\op{Com}C\ar@<1ex>[u]^{L:=?\otimes_{\tau}A}
}
\]
In particular, we have isomorphisms
\[(\ch A)(LN,M)\cong\H 0F_{\tau}(\HOM_{k}(N,M))\cong(\ch C)(N,RM)
\]
natural in $N$ and $M$.
\end{corollary}
\bigskip

\section{The bar derived category}\label{The bar derived category}

Let $A$ be a unital cdg algebra. When $C$ is the bar construction $BA$
of $A$, and the twisting cochain is the composition of the projection
with the shift to the right, $\tau_{A}:BA\arr{p_{1}}SA\arr{w}A$, we
will show that the conditions of the Corollary \ref{Factorizations
provided by adjunctions} are satisfied.

\subsection{The bar resolution of a cdg module}\label{The bar resolution of a cdg module}

In \cite[Lemme 2.2.1.9]{Lefevre03} it is essentially proved that given
a unital dg module $M$ over a unital dg algebra $A$, the adjunction
morphism
\[
\delta_{M}: LRM\ra M
\]
is a homotopy equivalence of complexes over $k$. This morphism is
usually called \emph{the bar resolution of $M$}.  We claim that in the
presence of curvature, the morphism $R(\delta_{M})$ is still a
homotopy equivalence. Indeed, if $sx$ and $sy$ are homogeneous
elements of $BA$, the morphism of graded $k$-modules
\[
g:(RLR)(M)\ra (LR)(M)\ko m\otimes sx\otimes a\otimes
sy\mapsto(-1)^{|m|+1}m\otimes(sx,sa,sy)\otimes 1_{A}
\]
induces a contracting homotopy for $\id-\eta_{RM}R(\delta_{M})$. That
is to say, $R(\delta_{M})$ has injective kernel.
We could say that $R(\delta_{M})$ hides a bar resolution of $M$. The precise statement is that $\delta_{M}: LRM\ra M$ is a cofibrant approximation in the model category $\cc A$, endowed with the model structure described in Theorem \ref{Factorizations provided by adjunctions}.

Now, from the identity $\delta_{LN}L(\eta_{N})=\id_{LN}$ we get
$R(\delta_{LN})RL(\eta_{N})=\id_{RLN}$, and since $R(\delta_{LN})$ is
a stable isomorphism, the same holds for $RL(\eta_{N})$. Therefore,
$L(\eta_{N})$ is a weak equivalence in $\cc A$. But it is a morphism
between fibrant-cofibrant objects, which implies that it is a homotopy
equivalence. Since left homotopy agrees with homotopy in the sense of
Frobenius categories (\emph{cf.}\ Theorem \ref{A recognition criterion in Frobenius categories}), we have that $L(\eta_{N})$ is a stable
isomorphism, \ie has injective cokernel.

Therefore, subsection \ref{Factorizations provided by adjunctions}
tells us that $\cc A$ and $\op{Com}BA$ admit certain model structures
making the bar/cobar adjunction into a Quillen equivalence.
We have proved the

\begin{theorem}
Let $A$ be a unital cdg algebra. Consider the bar/cobar adjunction
\[\xymatrix{\cc A\ar@<1ex>[d]^R \\
\op{Com}BA\ar@<1ex>[u]^L
}
\]
by using the twisting cochain $\tau_{A}$. Then
\begin{enumerate}[1)]
\item There is a model structure in $\cc A$ such that all the objects are fibrant, an object $M$ is cofibrant if and only if it is a direct summand of a module of the form $LN$ for some comodule $N$, and a morphism $f$ is a weak equivalence if and only if $Rf$ is a stable isomorphism.
\item There is a model structure in $\op{Com}BA$ such that all the objects are cofibrant, an object $N$ is fibrant if and only if it is a direct summand of a comodule of the form $RM$ for some module $M$, and a morphism $f$ is a weak equivalence if and only if $Lf$ is a stable isomorphism.
\item The bar/cobar adjunction is a Quillen equivalence for these model structures.
\end{enumerate}
\end{theorem}

\subsection{The bar derived category}\label{Definition of the bar derived category}

The homotopy category of $\cc A$ regarded with the model structure of
Theorem \ref{The bar resolution of a cdg module} is the
\emph{bar derived category} of $A$, denoted by $\cd_{bar}A$. It is the
localization of $\cc A$ with respect to the class of morphisms $f$
such that $R(f)$ is a stable isomorphism in $\op{Com}BA$, \ie such
that $\cone(Rf)=R\cone(f)$ is injective. Put differently, $\cd_{bar}A$
is the localization of $\cc A$ with respect to the class of morphisms
$f$ such that $\cone(f)$ is contractible as an $A_{[0,\infty[}$-module
over $A$.

Moreover, as we know from Theorem \ref{Factorizations provided by adjunctions}, the `acyclic modules' for the bar derived category, called \emph{bar acyclic modules}, are those $M$ such that
\[(\ch A)(LN,M)=0
\]
for every counital dg $BA$-comodule $N$. Hence, they are indeed those modules contractible as $A_{[0,\infty[}$-modules over $A$, but they admit (\emph{cf.}\ Corollary \ref{The adjunction between Frobenius categories}) an alternative description in terms of true acyclicity of certain complexes forming a class: the bar acyclic modules are those $M$ such that
\[\H 0F_{\tau}(\HOM_{k}(N,M))=0
\]
for every counital dg $BA$-comodule $N$. Accordingly, we use the term
\emph{bar quasi-isomorphisms} for 
the weak equivalences of $\cc A$, \ie
those morphisms $f$ such that $\cone(f)$ is bar acyclic, and the
\emph{bar closed} modules are those $M$ such that
\[(\ch A)(M,M')=0
\]
for every bar acyclic $M'$. Theorem \ref{Factorizations provided by
  adjunctions} states that $\cd_{bar}A$ is triangle equivalent to the
full subcategory $\ch_{p,bar}A$ of $\ch A$ formed by the bar closed
modules.
\bigskip

\section{The various derived categories of a dg algebra}\label{The various derived categories of a dg algebra}

Let $A$ be a unital dg $k$-algebra (regarded as a cdg algebra with vanishing curvature). \emph{A priori}, we have three derived categories associated to it, each of them made from its own notion of acyclicity: the classical, the relative and the bar derived category. But it turns out that these last two categories agree.

\subsection{The classical derived category}\label{The classical derived category}

A unital dg $A$-module $M$ is \emph{acyclic} if (when $k$ is trivially made into a counital dg $BA$-comodule) we have
\[(\ch A)((Lk)[n],M)\cong(\ch A)(A[n],M)\cong\H {-n}M=0
\]
for each $n\in\Z$, \ie if the complex $F_{\tau}(\HOM_{k}(k,M))\cong M$ is acyclic. A morphism $f$ is a \emph{quasi-isomorphism} if its cone is acyclic,
and a module $M$ is \emph{closed} if
\[(\ch A)(M,M')=0
\]
for every acyclic $M'$. As we know (\emph{cf.}\ Example \ref{Small object argument in Frobenius categories}), the \emph{derived category} of $A$, denoted by $\cd A$, is the localization of $\cc A$ with respect to the quasi-isomorphisms, and it is triangle equivalent to the full triangulated subcategory $\ch_{p}A$ of $\ch A$ formed by the closed modules.

\subsection{The relative derived category}\label{The relative derived category}

A unital dg $A$-module $M$ is \emph{relatively acyclic} if
\[(\ch A)((LK)[n],M)\cong(\ch A)((K\otimes_{\tau}A)[n],M)\cong\H {-n}F_{\tau}(\HOM_{k}(K,M))=0
\]
for each $n\in\Z$ and each $k$-module $K$ (trivially made into a counital dg $BA$-comodule). But $F_{\tau}(\HOM_{k}(K,M))$ is the complex $\HOM_{k}(K,M)$ with the differential induced by that of $M$. Hence, to be relatively acyclic amounts to being contractible as a complex over $k$. A morphism $f$ is a \emph{relative quasi-isomorphism} if its cone is relatively acyclic, and a module $M$ is \emph{relatively closed} if
\[(\ch A)(M',M)=0
\]
for every relatively acyclic $M'$. As we know \cite[Proposition 7.4]{Keller1998a}, the \emph{relative derived category} of $A$, denoted by $\cd_{rel}A$, is the localization of $\cc A$ with respect to the relative quasi-isomorphisms, and it is triangle equivalent to the full triangulated subcategory $\ch_{p,rel}A$ of $\ch A$ formed by the relatively closed modules.

But, in fact, the relative derived category of $A$ is its bar derived
category. Indeed, any bar acyclic dg module is certainly relatively
acyclic. On the other hand, let $M$ be a (unital) dg $A$-module which
is contractible as a complex over $k$, and let $h_{1}:M\ra M$ be the
contracting homotopy of that complex. By using Obstruction Theory (\emph{cf.}\
\cite[Chapitre B]{Lefevre03}), we can extend $h_{1}$ to a contracting
homotopy for $M$ regarded as an $A_{\infty}$-module over $A$, which
proves that $M$ is bar acyclic. Then, $\cd_{rel}A$ is the homotopy
category of a certain model category, and Theorem
\ref{Factorizations provided by adjunctions} proves the analogue of
\cite[Proposition 7.4]{Keller1998a} and tells us that $\ch_{p,rel}A$
is the smallest full triangulated subcategory of $\ch A$ containing
all the modules of the form $N\otimes_{\tau}A$, for any counital dg
$BA$-comodule $N$, and closed under arbitrary coproducts.

\subsection{Links between the derived categories}\label{Links between the derived categories and an alternative proof of Kenji's Theorem}

From the inclusion
\[\{(Lk)[n]\}_{n\in\Z}\subseteq\{LN\}_{n\in\Z\ko N\in\op{Com}BA}
\]
we have that
\[\{\text{relatively acyclic}\}\subseteq\{\text{acyclic}\},
\]
and thus $\ch_{p}A$ is a full triangulated subcategory of $\ch_{p,bar}A$. Therefore, we have a fully faithful triangulated functor
\[\cd A\ra\cd_{rel}A\ko M\mapsto \textbf{p}M
\]
where $\textbf{p}M$ is the closed resolution of $M$. It is well known that, if $k$ is a field, then the relatively acyclic unital dg $A$-modules are precisely the acyclic unital dg $A$-modules, and so we have $\cd A=\cd_{rel}A$. Then the following result should be viewed as an analogue of Kenji Lef{\`e}vre-Hasegawa's theorem \cite[Th\'{e}or\`{e}me 2.2.2.2]{Lefevre03} for unital dg algebras over an arbitrary commutative ring:

\begin{theorem}
Let $A$ be a unital dg $k$-algebra. Endow the category $\cc A$ of unital dg $A$-modules with its structure of model category whose homotopy category is the relative derived category $\cd_{rel}A$ (\emph{cf.}\ subsection \ref{The relative derived category}). Then there exists a model structure in $\op{Com}BA$ such that the bar/cobar constructions
\[\xymatrix{ \cc A\ar@<1ex>[d]^{R} \\
\op{Com}BA\ar@<1ex>[u]^{L}
}
\]
yield a Quillen equivalence.
\end{theorem}

Our proof of this theorem, being more conceptual than that of \cite{Lefevre03}, suggests that the classical derived category appears in the statement of Kenji Lef{\`e}vre-Hasegawa's theorem `incidentally', due to the fact that he is working over a semisimple base. In fact, the category that appears in its own right is the relative derived category.
\bigskip

\section{The relative derived category regarded from the $A_{\infty}$-theory}

The spirit of $A_{\infty}$-theory (at least over a field) is to replace quasi-isomorphisms by homotopy equivalences, up to increasing the amount of morphisms and/or the amount of objects. Key examples of that are:
\begin{enumerate}[1)]
\item The `Th\'{e}or\`{e}me des $A_{\infty}$-quasi-isomorphismes' \cite[Corollaire 1.3.1.3]{Lefevre03}, which states that if $k$ is a field, the category of dg $k$-algebras up to quasi-isomorphisms is equivalent to the category of dg $k$-algebras (with $A_{\infty}$-morphisms) up to homotopy equivalences of $A_{\infty}$-algebras.
\item The corresponding result for modules \cite[Proposition 2.4.1.1]{Lefevre03}, which states that if $k$ is a field and $A$ is an augmented dg $k$-algebra, then the derived category of $A$ is equivalent to the category of unital dg $A$-modules (with strictly unital morphisms of $A_{\infty}$-modules) up to strictly unital homotopy equivalences of  $A_{\infty}$-modules.
\end{enumerate}
In the proofs of these results the presence of a base field has been crucial. The aim of this section is to present an analog of the second one for an arbitrary commutative ring.

\subsection{Using arbitrary $A_{\infty}$-modules}

Let $A$ be an augmented dg $k$-algebra, and let $\tau:BA\ra A$ be the twisting cochain of section \ref{The bar derived category}. From the situation
\[\xymatrix{ \cc A\ar@<1ex>[d]^R\ar[rr]^{V\text{ (not full)}} && \op{Nod}_{\infty}A\ar[dll]^{\ \ \ \ \ \ R_{\infty}\text{ (fully faithful)}}\\
\op{Com}BA\ar@<1ex>[u]^L &&
}
\]
where $V$ takes a unital dg $A$-module $M$ to $M$ itself regarded as an $A_{\infty}$-module over $A$, we get an adjoint pair of functors
\[\xymatrix{ \cc A\ar@<1ex>[d]^{V} \\
\op{Nod}_{\infty}A\ar@<1ex>[u]^{LR_{\infty}}
}
\]
Notice that $\op{Nod}_{\infty}A$ inherits \emph{via} $R_{\infty}$ a structure of Frobenius category (use the structure of Frobenius category of $\op{Com}BA$ and the cone of a morphism between $A_{\infty}$-modules \cite[subsection 2.4.3]{Lefevre03}). Thus, we can consider the projective model structure on $\op{Nod}_{\infty}A$ (\emph{cf.}\ subsection \ref{Factorizations provided by adjunctions}) such that its associated homotopy category is the stable category $\ul{\op{Nod}}_{\infty}A$ of $\op{Nod}_{\infty}A$ up to homotopy equivalences of $A_{\infty}$-modules. Consider also $\cc A$ as a model category whose homotopy category is the relative derived category $\cd_{rel}A$ (\emph{cf.}\ section \ref{The relative derived category}).

\begin{theorem}
The adjunction
\[\xymatrix{\cc A\ar@<1ex>[d]^{V} \\
\op{Nod}_{\infty}A\ar@<1ex>[u]^{LR_{\infty}}
}
\]
is a Quillen equivalence. In particular, we have a triangle equivalence $\cd_{rel}A\simeq\ul{\op{Nod}}_{\infty}A$.
\end{theorem}
\begin{proof}
The only not straightforward step is to check that the unit of the adjunction $\eta_{M}\in(\op{Nod}_{\infty}A)(M,M\otimes_{\tau}BA\otimes_{\tau}A)$ is a weak equivalence in $\op{Nod}_{\infty}A$. Since $A$ is augmented, we have that $(\eta_{M})_{1}:M\ra M\otimes_{\tau}BA\otimes_{\tau}A$ is an inflation in the category $\cc k$ of complexes over $k$ (with the Frobenius structure given by the degreewise split short exact sequences). We can see that its cokernel is contractible in $\cc k$ with the contracting homotopy induced by the map
\[H(m\otimes sx\otimes a)=\begin{cases}(-1)^{|m|+|sx|}m\otimes(sx,sa)\otimes 1_{A} & \text{if }a\neq 1_{A}, \\ 0 & \text{if }a=1_{A}. \end{cases}
\]
where $sx$ is a homogeneous element of $BA$. Therefore (see \emph{e.g.} \cite[Lemme 4.18]{Cisinski2003}) $(\eta_{M})_{1}$ is a homotopy equivalence of complexes over $k$, that is to say, its cone $\cone((\eta_{M})_{1})$ is a contractible complex. But this is the underlying complex of the cone $\cone\eta_{M}$ of $\eta_{M}$, and by Obstruction Theory we conclude that $\cone\eta_{M}$ is a contractible $A_{\infty}$-module. In other words, $\eta_{M}$ is a weak equivalence in $\op{Nod}_{\infty}A$.
\end{proof}

Notice that the passage from $\cc A$ to $\op{Nod}_{\infty}A$ increases both the number of objects and of morphisms. However, it is not necessary to increase the number of objects, as shown in the following result:

\begin{corollary}
$\cd_{rel}A$ is triangle equivalent to the category of unital dg $A$-modules (with morphisms of $A_{\infty}$-modules) up to homotopy equivalences of $A_{\infty}$-modules.
\end{corollary}
\begin{proof}
Let $\cc$ be the full subcategory of $\op{Nod}_{\infty}A$ formed by the unital dg $A$-modules. It is an exact subcategory of $\op{Nod}_{\infty}A$ which inherits a structure of Frobenius category. Then, the inclusion $\cc\hookrightarrow\op{Nod}_{\infty}A$ induces a fully faithful functor
$\ul{\cc}\hookrightarrow\ul{\op{Nod}}_{\infty}A$ which is essentially surjective. Indeed, given an $A_{\infty}$-module $M$ over $A$, we know from the proof of the theorem above that the unit of the adjunction $(LR_{\infty},V)$,
\[\eta_{M}:M\ra M\otimes_{\tau}BA\otimes_{\tau}A,
\]
is a homotopy equivalence of $A_{\infty}$-modules, and so $M\cong M\otimes_{\tau}BA\otimes_{\tau}A$ in $\ul{\op{Nod}}_{\infty}A$. Therefore,
\[\cd_{rel}A\simeq\ul{\op{Nod}}_{\infty}A\simeq\ul{\cc}.
\]
\end{proof}

\subsection{Using strictly unital $A_{\infty}$-modules}

Let $A$ be an augmented dg $k$-algebra. If we want strictly unital morphisms and homotopies of $A_{\infty}$-modules to appear in the description of $\cd_{rel}A$ from the viewpoint of $A_{\infty}$-theory, then we have to use a different coalgebra. Namely, if $\ol{A}$ is the reduction of $A$, then we have to consider the bar construction $B\ol{A}$ of $\ol{A}$ instead of the bar construction $BA$ of $A$. The coalgebra $B\ol{A}$ is a counital dg coalgebra and the composition
\[\tau:B\ol{A}\arr{p_{1}}S\ol{A}\arr{w}\ol{A}\hookrightarrow A
\]
is a solution of the Maurer-Cartan equation of the dg algebra $\HOM_{k}(B\ol{A},A)$. Hence, we have an adjunction
\[\xymatrix{\cc A\ar@<1ex>[d]^{R} \\
\op{Com}B\ol{A}\ar@<1ex>[u]^L
}
\]
We can complete the picture with the functor $U:\cc A\ra\Mod_{\infty}A$ which takes a unital dg $A$-module $M$ to $M$ itself regarded as a strictly unital $A_{\infty}$-module over $A$, the equivalence $\ol{?}:\Mod_{\infty}A\arr{\sim}\op{Nod}_{\infty}A$ which takes a strictly unital $A_{\infty}$-module $M$ over $A$ to the $A_{\infty}$-module $M$ over $\ol{A}$ with restricted multiplications, and the bar construction $R_{\infty}:\op{Nod}_{\infty}\ol{A}\ra\op{Com}B\ol{A}$. Then we have
\[\xymatrix{\cc A\ar[rr]^{U\text{ (not full)}}\ar@<1ex>[d]^R && \Mod_{\infty}A\ar[r]_{\sim}^{\ol{?}} & \op{Nod}_{\infty}\ol{A}\ar[dlll]^{\ \ \ \ \ R_{\infty}\text{ (fully faithful)}} \\
\op{Com}B\ol{A}\ar@<1ex>[u]^L &&&
}
\]
If $\delta_{M}:LRM\ra M$ is the counit of the adjunction $(L,R)$, we can prove that $R(\delta_{M})$ is a homotopy equivalence of comodules by using almost the same contracting homotopy as the one used in subsection \ref{The bar resolution of a cdg module} (from which one can easily deduce the contracting homotopy -in the category of complexes over $k$- of \cite[Lemme 2.2.1.9]{Lefevre03}). Therefore, Theorem \ref{Factorizations provided by adjunctions} tells us that there exist certain model structures on $\cc A$ and $\op{Com}B\ol{A}$ making $(L,R)$ into a Quillen equivalence. It is easy to check that the homotopy category $\op{Ho}\cc A$ is again the relative derive category $\cd_{rel}A$. As before,
\[\xymatrix{\cc A\ar@<1ex>[d]^{U} \\
\Mod_{\infty}A\ar@<1ex>[u]^{L\circ R\circ\ol{?}}
}
\]
is an adjunction and $\Mod_{\infty}A$ is a Frobenius category with homotopies given by strictly unital homotopies of $A_{\infty}$-modules.

The proof of the following results are similar to those of the corresponding results above.

\begin{theorem}
The adjunction
\[\xymatrix{\cc A\ar@<1ex>[d]^{U} \\
\Mod_{\infty}A\ar@<1ex>[u]^{L\circ R_{\infty}\circ\ol{?}}
}
\]
is a Quillen equivalence. In particular, we have a triangle equivalence $\cd_{rel}A\simeq\ul{\Mod}_{\infty}A$.
\end{theorem}

\begin{corollary}
$\cd_{rel}A$ is triangle equivalent to the category of unital dg $A$-modules (with morphisms of strictly unital $A_{\infty}$-modules) up to strictly unital homotopy equivalences of $A_{\infty}$-modules.
\end{corollary}

\end{document}